\documentclass{amsart}
\usepackage{epsfig, eepic, epic}
\usepackage{color}
\usepackage{amssymb,amsmath, amsfonts}
\newtheorem{theorem}{Theorem}[section]

\newtheorem{corollary}[theorem]{Corollary}

\newtheorem{realizabilityconjecture}[theorem]{Realizability Conjecture}
\newtheorem{lemma}[theorem]{Lemma}
\newtheorem{example}[theorem]{Example}

\newtheorem{criterion}[theorem]{Criterion}
\newcommand{\PSbox}[3]{\mbox{\includegraphics{#1}\hspace{#2}\rule{0pt}{#3}}}

\newcommand{\chs}[2]{
\left(\begin{smallmatrix} #1 \\#2 
        \end{smallmatrix} \right)}
\def\l{\ell}
\def\wnot{w_{\rm o}}
\thicklines
\newcommand{\doubleci}{\put(0,0){\circle*{7}}\put(-1,0){\circle{15}}}
\newcommand{\ci}{\bullet}

\newcommand{\emt}{\mbox{ }}

\newcommand{\cinum}[1]{\put(5,5){\circle{15}} #1}

\newcommand{\rkj}{\mathrm{rk}_{j}}
\newcommand{\rk}{\mathrm{rk}}
\newcommand{\id}{\mathrm{id}}
\newcommand{\flags}{\mathcal{F}{l}_n}
\newcommand{\meet}{\wedge}
\newcommand{\join}{\vee}

\newcommand{\bigjoin}{\bigvee}
\newcommand{\given}{\,\, | \,\,}

\newcommand{\transverse}{T}

\newcommand{\Q}{\mathbb{Q}}
\newcommand{\proj}{\mathbb{P}}
\newlength{\cellsize} 
\newsavebox{\cell}
\sbox{\cell}{\begin{picture}(18,18)
\put(0,0){\line(1,0){18}}
\put(0,0){\line(0,1){18}}
\put(18,0){\line(0,1){18}}
\put(0,18){\line(1,0){18}}
\end{picture}}
\newcommand\cellify[1]{\def\thearg{#1}\def\nothing{}%
\ifx\thearg\nothing
\vrule width0pt height\cellsize depth0pt\else
\hbox to 0pt{\usebox{\cell} \hss}\fi%
\vbox to \cellsize{
\vss
\hbox to \cellsize{\hss$#1$\hss}
\vss}}
\newcommand\tableau[1]{\vtop{\let\\\cr
\baselineskip -16000pt \lineskiplimit 16000pt \lineskip 0pt
\ialign{&\cellify{##}\cr#1\crcr}}}
\begin{document}

\bibliographystyle{apalike}

\title[Intersections of Schubert varieties]{Intersections of Schubert 
varieties  and \\
other permutation array schemes}

\author{Sara Billey and Ravi Vakil}

\address{Department of Mathematics, University of Washington, Seattle, WA}
\email{billey@math.washington.edu}

\address{Department of Mathematics, Stanford University, Stanford, CA}
\email{vakil@math.stanford.edu}
%\urladdr{http://www.math-mit.edu/\mytilde sara/\medskip}
%{\it 
%\textup{2000} Mathematics Subject Classification:} 14M15.

%\subjclass[2000]{Primary 14M17; Secondary 14M15.}
\date{Wednesday, February 16, 2005.} 

\thanks{Supported by NSF grants DMS-9983797 and DMS-0238532
respectively.}

\begin{abstract}
  Using a blend of combinatorics and geometry, we give an algorithm
  for algebraically finding all flags in any zero-dimensional
  intersection of Schubert varieties with respect to three transverse
  flags, and more generally, any number of flags.  In particular, the
  number of flags in a triple intersection is also a structure
  constant for the cohomology ring of the flag manifold.  Our
  algorithm is based on solving a limited number of determinantal
  equations for each intersection (far fewer than the naive approach).
  These equations may be used to compute Galois and monodromy groups
  of intersections of Schubert varieties.  We are able to limit the
  number of equations by using the \textit{permutation arrays} of
  Eriksson and Linusson, and their permutation array varieties,
  introduced as generalizations of Schubert varieties.  We show that
  there exists a unique permutation array corresponding to each
  realizable Schubert problem and give a simple recurrence to compute
  the corresponding rank table, giving in particular a simple
  criterion for a Littlewood-Richardson coefficient 
  to be $0$.  We describe pathologies of Eriksson and
  Linusson's permutation array varieties (failure of existence,
  irreducibility, equidimensionality, and reducedness of equations),
  and define the more natural {\em permutation array schemes}.  In
  particular, we give several counterexamples to the Realizability
  Conjecture based on classical projective geometry.  Finally, we give
  examples where Galois/monodromy groups experimentally appear to be
  smaller than expected.
\end{abstract}

\maketitle

\tableofcontents

\section{Introduction}

A typical \textit{Schubert problem} asks how many lines in three-space
meet four generally chosen lines.  The answer, two, may be obtained by
computation in the cohomology ring of the Grassmannian variety
of two-dimensional planes in four-space.  Such questions were
considered by H.~Schubert in the nineteenth century.  During the past
century, the study of the Grassmannian has been generalized to the
flag manifold where one can ask analogous questions.

The flag manifold $\flags(K)$ parameterizes the complete flags
  $$F_{\ci}= \{ \{ 0 \} = F_0 \subset F_{1} \subset \cdots \subset
   F_{n}=K^{n}\}$$
where $F_{i}$ is a vector space of dimension $i$.  
({\em Unless otherwise noted, we will work over an arbitrary base field
  $K$.}  The reader, and Schubert, is welcome to assume $K=\mathbb{C}$
throughout.  For a general field, we should use the Chow ring rather
than the cohomology ring, but they agree for $K=\mathbb{C}$.  For this reason,
and in order not to frighten the reader, we will use the term
``cohomology'' throughout.)

A modern Schubert problem asks how many flags have relative position
$u,v,w$ with respect to three generally chosen fixed flags $X_{\ci}$,
$Y_{\ci}$ and $Z_{\ci}$.  The solution to this problem, due to Lascoux
and Sch\"{u}tzenberger \cite{LS1}, is to compute a product of Schubert
polynomials and expand in the Schubert polynomial basis.  The
coefficient indexed by $u,v,w$ is the solution.  This corresponds to a
computation in the cohomology (or Chow) ring of the flag variety.
(Caution: this solution is known to work only in characteristic $0$,
due to the failure of the Kleiman-Bertini theorem, cf.\ \cite[Sect.\ 
2]{VakilB}.)  The quest for a combinatorial rule for expanding these
products is a long-standing open problem.

The main goal of this paper is to describe a method for directly
identifying all flags in $X_{u}(F_{\ci}) \cap X_{v}(G_{\ci}) \cap
X_{w}(H_{\ci}) $ when $\l (u) + \l (v) + \l (w) = \chs{n}{2}$, thereby
computing $c_{u,v,w}$ explicitly.  This method extends to Schubert
problems with more than three flags.  To do this, we use the
permutation arrays defined by Eriksson and Linusson to obtain a table
of intersection dimensions.  These permutation arrays are closely
related to the checker boards used in \cite{VakilA,VakilB}.  Eriksson
and Linusson introduced permutation array varieties as natural
generalizations of Schubert varieties to an arbitrary number of flags.
We show that they may be badly behaved. For example, their equations
are not reduced, so we argue that the ``correct'' generalization of
Schubert varieties are permutation array {\em schemes}.  We describe
pathologies of these varieties/schemes, and show that they are
not irreducible nor even equidimensional in general, making a
generalization of the Bruhat order problematic.  We also give
counterexamples to the Realizability Conjecture~\ref{c:realizability}.
Returning to the task at hand, we use the data from the permutation
array to identify and solve a collection of determinantal equations
for the permutation array schemes, allowing us to solve Schubert
problems explicitly and effectively, for example allowing us to
compute Galois/monodromy groups.

The outline of the paper is as follows.  In Section~\ref{s:flags}, we
review Schubert varieties and the flag manifold. 
In Section~\ref{s:perm-arrays}, we
review the construction of permutation arrays and the
Eriksson-Linusson algorithm for generating all such arrays.  
In Section~\ref{s:realizability}, we describe 
permutation varieties and their pathologies, and explain why their
correct definition is as schemes.
In
Section~\ref{s:intersections}, we describe how to use permutation
arrays to solve Schubert problems and give equations for certain
intersections of Schubert varieties.    In Section~\ref{s:triple},
we give two examples 
of an algorithm for computing triple intersections of Schubert
varieties and thereby computing the cup product in the cohomology ring
of the flag manifold.
The equations we give also allow
us to compute Galois and monodromy groups for intersections of
Schubert varieties; we describe this application in
Section~\ref{s:monodromy}.  Our computations lead to examples where
the Galois/monodromy group is smaller than expected.

\section{The Flag Manifold and Schubert varieties}\label{s:flags}

In this section we briefly review the notation and basic concepts for
flag manifolds and Schubert varieties.  We refer the reader to one of
the following books for further background information:
\cite{Fulton-book,M2,manivel-book,GL-book,kumar-book}.  

As described earlier, the flag manifold
$\flags=\flags(K)$ parametrizes the  complete flags
$$F_{\ci}= \{ \{ 0 \} = F_0 \subset F_{1} \subset \dots \subset
F_{n}=K^{n}\}$$ where $F_{i}$ is a vector space of dimension $i$ over
the field $K$.  $\flags$ is a smooth projective variety of dimension
$\chs{n}{2}$. A complete flag is determined by an ordered basis
$(f_{1},\dots, f_{n})$ for $K^{n}$ by taking $F_{i} =
\operatorname{span}(f_{1},\dots, f_{i})$.

Two flags $[F_{\ci}], [G_{\ci}] \in \flags$ are {\em in relative
position $w \in S_{n}$} when
\[
\dim(F_{i} \cap G_{j}) = \mathrm{rank}\, w[i,j] \hspace{.2in} 
\text{ for all } \ 1\leq i,j\leq n
\]
where $w[i,j]$ is the principal submatrix of the permutation matrix
for $w$ with lower right hand corner in position $(i,j)$.  
We use the
notation 
\[
\mathrm{pos}(F_{\bullet},G_{\bullet}) =w.  
\]
Warning: in order to use the typical meaning for a principal submatrix we
are using  a nonstandard labeling of a permutation matrix.
The permutation matrix we associate to $w$ has a $1$ in the $w(i)$th
row of column $n-i+1$ for $1\leq i\leq n$.  For example, the matrix
associated to $w=(5,3,1,2,4)$ is
\[
\begin{array}{ccccc}
0&0&1&0&0\\
0&1&0&0&0\\
0&0&0&1&0\\
1&0&0&0&0\\
0&0&0&0&1.\\
\end{array}
\]
If $\mathrm{pos}(F_{\bullet},G_{\bullet}) =(5,3,1,2,4)$ then
$\dim(F_{2} \cap G_{3})=2$ and $\dim(F_{3} \cap
G_{2})=1$.

Define a \textit{Schubert cell} with respect to a fixed flag $F_{\ci}$
in $\flags$ to be
%\begin{align}\label{e:schubert.cell}
\begin{eqnarray*}
X^{o}_{w}(F_\ci) &=& \{G_{\ci} \given F_{\ci} \text{ and } G_{\ci}
\text{ have relative position } w \}\\
&=& \{G_{\ci} \given \dim (F_i \cap  G_j ) =
\mathrm{rk}\, w[i,j] \}.
\end{eqnarray*}
Using our labeling of a permutation matrix, the codimension of
$X^{o}_{w}$ is equal to the length of $w$ (the number of inversions in $w$), denoted $\l (w)$.  In fact, $X^{o}_{w}$ is
isomorphic to the affine space $K^{\chs{n}{2}-\l(w)}$.  We say the flags $F_{\ci}$ and $G_{\ci}$ are
in \textit{transverse position} if $G_{\ci} \in X_{\id}(F_{\ci})$.  A
randomly chosen flag will be transverse to any fixed flag $F_{\ci}$
with probability 1 (using any reasonable measure).

The \textit{Schubert variety} $X_{w}(F_{\ci})$ is the closure of
$X^{o}_{w}(F_{\ci})$ in $\flags$.  Schubert
varieties may also be written in terms of rank conditions:
%\begin{equation}\label{}
$$
X_{w}(F_{\ci}) = \{G_{\ci} \given \dim (F_i \cap  G_j ) \geq 
\mathrm{rk}\, w[i,j] \}.
$$
%\end{equation}
If the flags $F_{\ci}$ and $G_{\ci}$ are determined by ordered bases
for $K^{n}$ then these inequalities can be rephrased as
determinantal equations on the coefficients in the bases
\cite[10.5, Ex.~10, 11]{Fulton-book}.   Of course this allows
one in theory to solve all Schubert problems, but the number and
complexity of the equations grows quickly to make this prohibitive 
in any reasonable case.

%Ehresmann\cite{} and Chevalley\cite{CHEV} first studied the poset
%on permutations defined by inclusions of Schubert varieties which is
%now known as Bruhat order,
%\[
%X_{w}(F_{\ci}) = \bigcup_{v \leq w} X_{v}(F_{\ci}).
%\]
%The Bruhat order on permutations can also be defined by taking the
%transitive closure of the relation $v<vt_{ij}$ if $\l(v)<l(vt_{ij})$.

The cohomology (or Chow) ring of $\flags$ 
is isomorphic to 
$\mathbb{Z}[x_{1},\dots,x_{n}]/\langle e_{1},
e_{2}, \dots, e_{n} \rangle $ where $e_{i}$ is the $i$th elementary
symmetric function on $x_{1},\dots, x_{n}$
\cite[10.2, B.3]{Fulton-book}.   The cycles $[X_{u}]$
corresponding to Schubert varieties form a $\mathbb{Z}$-basis for the ring.
The product is defined by 
$$
%\begin{equation}\label{e:product}
[X_{u}] \cdot [X_{v}] = [X_{u}(F_{\ci}) \cap X_{v}(G_{\ci})] 
%\end{equation}
$$
where $F_{\ci}$ and $G_{\ci}$ are in {transverse position}.
Speaking somewhat informally, $X_{u}(F_{\ci}) \cap X_{v}(G_{\ci})$ can be decomposed into
irreducible components which are again %$GL_n$-translates of 
translates of 
Schubert varieties.  Therefore the expansion
\begin{equation}\label{e:coefs} % used
[X_{u}] \cdot [X_{v}] = \sum_{\l (w) = \l (u) + \l (v)} c_{u,v}^{w} [X_{w}]
\end{equation}
automatically has nonnegative integer coefficients.

A simpler geometric interpretation of the coefficients $c_{u,v}^{w}$
may be given in terms of triple intersections
\cite[10.2]{Fulton-book}.  There exists a perfect pairing on
$H^{*}(\flags)$ such that
\begin{equation}\label{e:compliments}  % used
[X_{w}]\cdot [X_{y}] = 
\begin{cases}
[X_{\wnot}]  &  y=\wnot w\\
0 &      y \neq \wnot w, \; \l (y) = \chs{n}{2} -\l (w).
\end{cases}
\end{equation}
Here $\wnot = (n,n-1,\dots, 1)$ is the longest permutation in
$S_{n}$, of  length  $\chs{n}{2}=\dim (\flags)$, and
$[X_{\wnot}]$ is the class of a point. 
% $X_{\wnot}$ is the unique Schubert variety of dimension zero.  
Combining equations
\eqref{e:coefs} and \eqref{e:compliments} we have
$$
%\begin{equation}\label{e:triple}
[X_{u}] \cdot [X_{v}]\cdot [X_{\wnot w}] = c_{u,v}^{w} [X_{\wnot}].  
%\end{equation}
$$
In characteristic $0$, $c_{u,v}^{w}$ counts the number of points
$[E_{\ci}] \in \flags$ in the variety
\begin{equation}\label{e:triple-intersection} % used
X_{u}(F_{\ci}) \cap X_{v}(G_{\ci}) \cap   X_{\wnot w}(H_{\ci}) 
\end{equation}
when $\l(u)+\l(v)+\l(\wnot w)=\chs{n}{2}$ and $F_{\ci}, G_{\ci},
H_{\ci}$ are three generally chosen flags. 
Note, it is not sufficient to assume the three flags are pairwise
transverse in order to get the expected number of points in the
intersection.  There can be additional dependencies among the
subspaces of the form $F_{i} \cap G_{j} \cap H_{k}$.

The main goal of this article is to describe a method to find all
flags in a general $d$-fold intersection of Schubert varieties when
the intersection is zero-dimensional.
%  For our purposes, those
%configurations of flags with the expected  number of solutions
%will be general enough.  (The number of solutions will always be
%infinite, or no greater than the expected number. 
%The expected number is
%achieved on a dense open subset of $\flags^3$.)  
Enumerating the flags
found explicitly in a triple intersection would give the numbers
$c_{u,v}^{w}$.  We will use the permutation arrays defined in the next
section to identify a different set of equations defining the
intersections of Schubert varieties which are easier to solve.

\section{Permutation arrays}\label{s:perm-arrays}

In \cite{ELcombinatorial,ELdecomposition}, Eriksson and Linusson
develop a $d$-dimensional analog of a permutation matrix.  One way to
generalize permutation matrices is to consider all $d$-dimensional
arrays of $0$'s and $1$'s with a single $1$ in each hyperplane.
They claim that a better way is to consider a permutation matrix to be
a two-dimensional array of 0's and 1's such that the rank of any
principal minor is equal to the number of occupied rows in that
submatrix or equivalently equal to the number of occupied columns in
that submatrix.  The locations of the 1's in a permutation matrix will
be the elements in the corresponding permutation array.  We will
summarize their work here and refer the reader to their well-written
paper for further details.

Let $P=\{(x_{1}, \dots, x_{d})\}$ be any collection of points in
$[n]^{d}:=\{1,2,\dots, n\}^{d}$.  We will think of these points as
the locations of dots in an $[n]^{d}$-\textit{dot array}.  Define a
partial order on $[n]^{d}$ by
$$x=(x_{1},\dots, x_{d}) \preceq y=(y_{1}, \dots, y_{d}),$$
read ``$x$ is \textit{dominated} by $y$'', if $x_{i} \leq y_{i}$ for all
$1\leq i \leq d$.  This poset is a lattice with meet and join
operation defined by

\begin{eqnarray*}
x \join y &=& z \quad \quad \text{ if  } z_{i} = \mathrm{max} (x_{i},
y_{i}) \text{ for all $i$}
\\
x \meet y &=& z \quad \quad \text{ if } z_{i} = \mathrm{min} (x_{i}, y_{i})
\text{ for all $i$.}
\end{eqnarray*}

These operations extend to any set of points $R$ by taking $\bigvee
R=z$ where $z_{i}$ is the the maximum value in coordinate $i$ over the
whole set, and similarly for $\bigwedge R$.

Let $P[y]=\{x \in P \given x \preceq y \}$ be the \textit{principal subarray}
of $P$ containing all points of $P$ which are dominated by $y$.  Define 
\[
\rkj P= \# \{1\leq k\leq n \given \text{there exists } x \in P \text{ with } x_{j}=k \}.
\]
$P$ is \textit{rankable} of \textit{rank} $r$ if $\rkj P =r$ for all
$1\leq j \leq d$.  $P$ is \textit{totally rankable} if every principal
subarray of $P$ is rankable.  

For example, with $n=4$, $d=3$ the following example is a totally
rankable dot array: $\{(3,4,1), (4,2,2), (1,4,3), (3,3,3), (2,3,4),
(3,2,4),(4,1,4)\}$.  We picture this as four 2-dimensional slices, where
the first one is ``slice $1$'' and the last is ``slice $4$'':
$$
\tableau{\emt & \emt & \emt & \emt \\
         \emt & \emt & \emt & \emt \\
         \emt & \emt & \emt & \ci \\
         \emt & \emt & \emt & \emt }
\hspace{.2in}
\tableau{\emt & \emt & \emt & \emt \\
         \emt & \emt & \emt & \emt \\
         \emt & \emt & \emt & \emt \\
         \emt & \ci & \emt & \emt }
\hspace{.2in}
\tableau{\emt & \emt & \emt & \ci \\
         \emt & \emt & \emt & \emt \\
         \emt & \emt & \ci & \emt \\
         \emt & \emt & \emt & \emt }
\hspace{.2in}
\tableau{\emt & \emt & \emt & \emt \\
         \emt & \emt & \ci      & \emt \\
         \emt & \ci      & \emt & \emt \\
         \ci      & \emt & \emt & \emt }
$$
Thus $(3,4,1)$ corresponds to the dot in the first slice on the left.

The array $\{(3,4,1), (4,2,2), (1,4,3)\}$ is not rankable since it has
only two distinct values appearing in the second index and three in
the first and third.  

Many pairs $P, P'$ of totally rankable dot arrays are \textit{rank
equivalent}, i.e.\ $\rkj P[x] =\rkj P'[x]$, for all $x$ and $j$.
However, among all rank equivalent dot arrays there is a unique one
with a minimal number of dots \cite[Prop.\ 4.1]{ELcombinatorial}.  In
order to characterize the minimal totally rankable dot arrays, we give
the following two definitions.  We say a position $x$ is
\textit{redundant} in $P$ if there exists a collection of points $R
\subset P$ such that $x = \bigvee R$, $\# R >1$, and every $y \in R$
has at least one $y_{i} =x_{i}$.  We say a position $x$ is
\textit{covered} by dots in $P$ if $x$ is redundant for some $R
\subset P$, $x \notin R$, and for each $1 \leq j\leq d$ there exists
some $y \in R$ such that $y_{j}<x_{j}$.  We show in
Lemma~\ref{l:covered-dots} that it suffices to check only subsets $R$
of size at most $d$ when determining if a position is redundant or
covered.  

\begin{theorem}\label{t:EL}\cite[Theorem 3.2]{ELdecomposition} Let $P$ be a dot array.  The
following are equivalent:
\begin{enumerate}
\item $P$ is totally rankable.
\item Every two dimensional projection of every principal subarray is
totally rankable.
\item Every redundant position is covered by dots in $P$.
\item If there exist dots in $P$ in positions $y$ and $z$ and integers
$i,j$ such that $y_{i}<z_{i}$ and $y_{j}=z_{j}$, then there exists a
dot in some position  
$x \preceq (y \join z)$ such that
$x_{i} = z_{i}$ and $x_{j} < z_{j}$.
\end{enumerate}
\end{theorem}

Define a \textit{permutation array} in $[n]^{d}$ to be a totally
rankable dot array of rank $n$ with no redundant dots (or equivalently,
no covered dots).  The permutation arrays are the unique
representatives of each rank equivalence class of totally rankable dot
arrays with no redundant dots.  These arrays are Eriksson and
Linussons' analogs of permutation matrices.

The definition of permutation arrays was motivated because they
include the possible relative configurations of flags:
\begin{theorem}  \label{t:EL-strata} 
\cite[Thm.\ 3.1]{ELdecomposition} 
Given flags $E_{\ci}^{1}, E_{\ci}^{2}, \dots , E_{\ci}^{d}$, there
exists an $[n]^{d}$-permutation array $P$ describing the table of all
intersection dimensions as follows.  For each $x \in [n]^{d}$,
\begin{equation}\label{e:EL-strata}
\rk(P[x]) = \dim \left(E_{x_{1}}^{1}\cap E_{x_{2}}^{2}\cap
\dots \cap E_{x_{d}}^{d} \right).
\end{equation}
\end{theorem}

A special case is the permutation array corresponding to $n$ 
generally chosen flags, which we denote the {\em transverse permutation 
array}
$$
\transverse_{n,d} = \left\{(x_{1},\dots , x_{d}) \in [n]^{d}\given \sum x_{i}
=(d-1)n+1 \right\},$$
which corresponds to 
$$
\rk(\transverse_{n,d} [x]) =
\max \left( 0,  n - \sum_{i=1}^d (n-x_i) \right).
$$

Eriksson and Linusson give an algorithm for producing all
permutation arrays in $[n]^{d}$ recursively from the permutation
arrays in $[n]^{d-1}$.  We review their algorithm, as this is key to
our algorithm for intersecting Schubert varieties.  

Let $A$ be any antichain of dots in $P$ under the dominance order.
Let $C(A)$ be the set of positions  covered by dots in 
$A$.  Define the \textit{downsizing} operator 
$D(A,P)$
with respect to $A$
on $P$  to be the result of the following process.
\begin{enumerate}
\item Set $Q_{1} = P \setminus A$. 
\item Set $Q_{2} = Q_{1} \cup C(A)$.
\item Set $D(A,P) = Q_{2} \setminus R(Q_{2})$ where $R(Q)$ is the set
of redundant positions of $Q$.
\end{enumerate}

The downsizing is \textit{successful} if the resulting array is
totally rankable of rank $\rk(P)-1$.

\begin{theorem}[The EL-Algorithm, {\cite[Sect.\ 2.3]{ELdecomposition}}]
\label{t:EL-algorithm}Every 
permutation array in $[n]^{d}$ can be obtained uniquely in the
following way.
\begin{enumerate}
\item Choose a permutation array $P_{n}$ in $[n]^{d-1}$.
\item For each $n\geq i>1$, choose an antichain $A_{i}$ of dots in
$P_{i}$ such that the downsizing $D(A_{i},P_{i})$ is successful.  Set
$P_{i-1}=D(A_{i}, P_{i})$.
\item Set $A_{1}=P_{1}$.  
\item Set $P = \{(x_{1},\dots, x_{d-1},i) \given
(x_{1},\dots,x_{d-1}) \in A_{i} \}$.  
\end{enumerate}
\end{theorem}

For example, starting with the $2$-dimensional array $\{(1,4),
(2,3),(3,1), (4,2) \}$ corresponding to the permutation $w=(1,2,4,3)$,
we run through the algorithm as follows.  (In the figure, dots
correspond to elements in $P$ and circled dots correspond to elements in $A$.)
\[
\begin{array}{ll}
P_{4}=\{(1,4), (2,3),(3,1), (4,2) \} &  A_{4} = \{(1,4), (2,3) \}\\
P_{3} = \{(2,4),(3,1), (4,2) \} &       A_{3}= \{(3,1) \}\\
P_{2} = \{(2,4), (4,2) \} &     A_{2}= \{(2,4), (4,2) \}\\
P_{1} = \{(4,4)\} &     A_{1}= \{(4,4)\} 
\end{array}
\]

$$
\tableau{\emt & \emt & \emt & \emt \\
         \emt & \emt & \emt & \emt \\
         \emt & \emt & \emt & \emt \\
         \emt & \emt & \emt & \doubleci }
\hspace{.2in}
\tableau{\emt & \emt & \emt & \emt \\
         \emt & \emt & \emt & \doubleci \\
         \emt & \emt & \emt & \emt \\
         \emt & \doubleci & \emt & \emt }
\hspace{.2in}
\tableau{\emt & \emt & \emt & \emt \\
         \emt & \emt & \emt & \ci \\
         \doubleci & \emt & \emt & \emt \\
         \emt & \ci & \emt & \emt }
\hspace{.2in}
\tableau{\emt & \emt & \emt & \doubleci \\
         \emt & \emt & \doubleci & \emt \\
         \ci & \emt & \emt & \emt \\
         \emt & \ci & \emt & \emt }
$$
This produces the 3-dimensional array 
\[
P= \{(4,4,1), (2,4,2), (4,2,2),
(3,1,3), (1,4,4), (2,3,4) \}.
\]

We prefer to display 3-dimensional dot-arrays as 2-dimensional
number-arrays as in \cite{ELdecomposition, VakilA} where a square $(i,j)$ contains
the number $k$ if $(i,j,k) \in P$.  Note that there is at most one number
in any square if the number-array represents a permutation array.  The
previous example is represented by 
$$
\tableau{\emt & \emt & \emt & 4 \\
         \emt & \emt & 4 & 2 \\
         3 & \emt & \emt & \emt \\
         \emt & 2 & \emt & 1 }.
$$
\vspace{.3in}
%The checkerboard notation makes it possible to represent 4-dimensional
%arrays as $n$ checkerboards.  

%Example:

%\begin{verbatim}
%    ....    ....    ...3    ...3
%    ....    ....    ....    ..32
%    ...3    ...3    ..31    .321
%    ....    .3.0    .310    3210
%\end{verbatim}

\begin{corollary}\label{c:checkerboards}
In Theorem~\ref{t:EL-algorithm}, each $P_{i}$ is an
$[n]^{d-1}$-permutation array of rank $i$.  Furthermore, if $P$
determines the rank table for flags $E^{1}_{\ci},\dots , E^{d}_{\ci}$,
then $P_{i}$ determines the rank table for $E^{1}_{\ci},\dots ,
E^{d-1}_{\ci}$ intersecting the vector space $E_{i}^{d}$, i.e.\
\[
\mathrm{rk}\left(P_{i}[x] \right) = \dim\left(E_{x_{1}}^{1}\cap
E_{x_{2}}^{2}\cap \cdots \cap E_{x_{d-1}}^{d-1} \cap E_{i}^{d} \right).
\]
\end{corollary}

\begin{proof}
$P_{i}$ is the permutation array obtained from the
projection
\[
\{(x_{1},\dots, x_{d}) \given (x_{1},\dots , x_{d}, x_{d+1}) \in P \text{
and } x_{d+1} \leq i \}
\]
by removing all repeated or covered elements.
\end{proof}

\bigskip

We finish this section with a substantial improvement on the speed to
the Eriksson-Linusson algorithm.  In Step~2 of
Theorem~\ref{t:EL-algorithm}, one must find all positions covered by a
subset of points in the antichain $A_{i}$.  This appears to require on
the order of $2^{|A_{i}|}$ computations.  However, here we show that
subsets of size at most $d$ are sufficient.

\begin{lemma}\label{l:covered-dots}
A position $x \in [n]^{d}$ is covered (or equivalently, redundant) in
a permutation array $P$ if and only if there exists a subset $S$ with
$|S|\leq d$ which cover $x$.
\end{lemma}

\begin{proof}
Assume $x$ is covered by a set $Y=\{y^{1}, y^{2}, \dots, y^{k} \}$ for
$k>d$.  That is,
\begin{itemize}
\item For each position $1 \leq j\leq d$, there exists a $y^{i}$ such
that $y^{i}_{j} < x_{j}$ and there exists a $y^{l}$ such that
$y^{l}_{j} = x_{j}$.
\item For each $y^{i} \in Y$, there exists a $j$ such that $y^{i}_{j} <
x_{j}$ and there exists an $l$ such that $y^{i}_{l} = x_{l}$.
\end{itemize}
Consider a complete bipartite graph with left vertices labeled by
$Y$ and right vertices labeled by $\{x_{1}, \dots, x_{d} \}$.  Color
the edge from $y^{i}$ to $x_{j}$ red if $y^{i}_{j} =x_{j}$, and blue if
$y^{i}_{j} <x_{j}$.  Since $x = \bigjoin Y$, $y^{i}_{j} >x_{j}$ is not
possible.  This is a complete bipartite graph such that each vertex
meets at least one red and one blue edge, and conversely
any such complete bipartite graph with left vertices chosen from $P$
and right vertices $\{x_{1}, \dots , x_{d} \}$ corresponds to a
covering of $x$.

We can easily bound the minimum size of a covering set for $x$ to be
at most $d+1$ as follows.  Choose one red and one blue edge adjacent
to $x_{1}$.  Let $S$ be the left end-points of these two edges.
Vertex $x_{2}$ is connected to both elements of $S$ in the complete
bipartite graph.  If the edges connecting $x_{2}$ to $S$ are different
colors, proceed to $x_{3}$.  If the edges agree in color, choose one
additional edge of a different color adjacent to $x_{2}$.  Add its
left endpoint to $S$.  Continuing in this way for $x_{3},\dots,
x_{d}$, we have $|S|\leq d+1$ and that $x$ is covered by $S$.

Given a covering set $S$ of size $d+1$, we now find a subset of size
$d$ which covers $x$.  Say $x_{i_{1}}, x_{i_{2}}, \dots, x_{i_{k}}$
are all the right vertices which are adjacent to a unique edge of
either color.  Let $T$ be the left endpoints of all of these edges;
these are necessary in any covering subset.  Choose one vertex in $Y
\setminus T$, say $\tilde{y}$.  Each remaining $x_{j}$ has at least
two edges of each color, so we can choose one of each color which is
not adjacent to $\tilde{y}$.  The induced subgraph on $(S \setminus
\{\tilde{y} \}, \{x_{1},\dots , x_{d} \})$ is again a complete
bipartite graph where every vertex is adjacent to at least one red and
one blue edge, hence $S \setminus \{\tilde{y} \}$ covers $x$.
\end{proof}

\section{Permutation array varieties (or schemes) and their pathologies}
\label{s:realizability}

In analogy with Schubert cells, for any $[n]^d$-permutation array
$P$, Eriksson and Linusson define the {\em permutation array variety}
$X^o_P$ to be the subset of $\flags^d = \left\{
\left( E^1_{\ci}, \dots, E^d_{\ci} \right) \right\}$ in ``relative
position $P$'' \cite[\S 1.2.2]{ELdecomposition}.  We will soon see why $X^o_P$ is a locally closed
subvariety of $\flags^d$; this will reinforce the idea that the correct notion is of
a permutation array {\em scheme}.
These varieties/schemes  will give a convenient way to 
manage the equations of intersections of Schubert varieties.

Based on many examples, Eriksson and Linusson conjectured the following.
\begin{realizabilityconjecture}[{\cite[Conj.~3.2]{ELdecomposition}}]
Every permutation array can be realized by flags.  Equivalently,
every $X^o_P$ is nonempty. \label{c:realizability}
\end{realizabilityconjecture}

This question is motivated by more than curiosity.  A fundamental
question is: {\em what are the possible relative configurations of $d$
  flags?}  In other words: {\em what intersection dimension tables are
  possible?}  For $d=2$, the answer leads to the theory of Schubert
varieties.  By Theorem~\ref{t:EL-strata}, each achievable intersection
dimension table yields a permutation array, and the permutation arrays
may be enumerated by Theorem~\ref{t:EL-algorithm}.  The Realizability
Conjecture then says that we have fully answered this fundamental
question.  Failure of realizability would imply that we still have a
poor understanding of how flags can meet.

The Realizability Conjecture is true for $d=1, 2, 3$.  For $d=1$, the only
permutation array variety is the flag variety.  For $d=2$, the
permutation array varieties are the 
``generalized'' Schubert cells (where the reference flag may vary).  The
case $d=3$ follows from \cite{SSV} (as described in \cite[\S
3.2]{ELdecomposition}), see also \cite[\S 4.8]{VakilA}.  The case $n
\leq 2$ is fairly clear, involving only one-dimensional subspaces of a
two-dimensional vector space (or projectively, points on $\proj^1$),
cf.\ \cite[Lemma~4.3]{ELdecomposition}.  Nonetheless, the conjecture
is false, and we give examples below which show the bounds $d \leq 3$
and $n \leq 2$ are maximal for such a realizability statement.  We
found it interesting that the combinatorics of permutation arrays
prevent some naive attempts at counterexamples from working; somehow,
permutation arrays see some subtle linear algebraic information, but
not all.

{\bf Fiber permutation array varieties.}  If $P$ is an $[n]^{d+1}$
permutation array, then there is a natural morphism $X^o_P \rightarrow
\flags^d$ corresponding to ``forgetting the last flag''.  We call the
fiber over a point $(E^1_{\ci}, \dots, E^d_{\ci})$ a {\em fiber
permutation array variety}, and denote it $X^o_P(E^1_{\ci}, \dots,
E^d_{\ci})$.  If the flags $E^1_{\ci}, \dots, E^d_{\ci}$ are chosen
generally, we call the fiber permutation array variety a {\em
generic fiber permutation array variety}.  Note that a generic fiber
permutation array variety is empty unless the projection of the
permutation array to the ``bottom hyperplane of $P$'' is the transverse
permutation array $\transverse_{n,d}$, as this projection describes
the relative positions of the first $d$ flags.

The Schubert cells $X_w^o(E^1_\ci)$ are fiber
permutation array varieties, with $d=2$.  Also, any
intersection of Schubert cells $$X_{w_1}(E^1_\ci) \cap
X_{w_2}(E^2_\ci) \cap \cdots \cap X_{w_d}(E^d_\ci)$$ is a
disjoint union of fiber permutation array varieties, and if the
$E^i_\ci$ are generally chosen, the intersection is a disjoint union
of generic fiber permutation array varieties.

Permutation array varieties were introduced partially for this reason, to
study intersections of Schubert varieties, and indeed that is the 
point of this paper.  It was hoped that
they would in general be tractable and well-behaved (cf.\ the Realizability
Conjecture~\ref{c:realizability}), but sadly this is not the case.  The remainder
of this section is devoted to their pathologies, and is independent
of the rest of the paper.

{\bf Permutation array schemes.}  We first observe that the more
natural algebro-geometric definition is of {\em permutation array
schemes}: the set of $d$-tuples of flags in configuration $P$ comes
with a natural scheme structure, and it would be naive to expect that
the resulting schemes are reduced.  In other words, the ``correct''
definition of $X^o_P$ will contain infinitesimal information not
present in the varieties.  More precisely, the $X^o_P$ defined above
may be defined scheme-theoretically by the equations
\eqref{e:EL-strata}, and these equations will {\em not} in general be
all the equations cutting out the {\em set} $X^o_P$ (see the ``Further
Pathologies'' discussion below).  Those readers preferring to avoid
the notion of schemes may ignore this definition; other readers should
re-define $X^o_P$ to be the scheme cut out by equations
\eqref{e:EL-strata}, which is a locally closed subscheme of
$\flags^d$.
%locally closed =  open subset of its closure

We now give a series of counterexamples to the Realizability
Conjecture~\ref{c:realizability}.

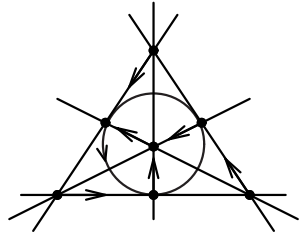
\begin{figure}[ht]
\begin{center}
\setlength{\unitlength}{0.00083333in}
\begingroup\makeatletter\ifx\SetFigFont\undefined%
\gdef\SetFigFont#1#2#3#4#5{%
  \reset@font\fontsize{#1}{#2pt}%
  \fontfamily{#3}\fontseries{#4}\fontshape{#5}%
  \selectfont}%
\fi\endgroup%
{\renewcommand{\dashlinestretch}{30}
\begin{picture}(1824,1464)(0,-10)
\put(912,554){\ellipse{634}{634}}
\put(912,1137){\blacken\ellipse{50}{50}}
\put(912,1137){\ellipse{50}{50}}
\put(1212,687){\blacken\ellipse{50}{50}}
\put(1212,687){\ellipse{50}{50}}
\put(1512,237){\blacken\ellipse{50}{50}}
\put(1512,237){\ellipse{50}{50}}
\put(912,237){\blacken\ellipse{50}{50}}
\put(912,237){\ellipse{50}{50}}
\put(312,237){\blacken\ellipse{50}{50}}
\put(312,237){\ellipse{50}{50}}
\put(612,687){\blacken\ellipse{50}{50}}
\put(612,687){\ellipse{50}{50}}
\put(912,537){\blacken\ellipse{50}{50}}
\put(912,537){\ellipse{50}{50}}
\path(162,12)(1062,1362)
\path(762,1362)(1662,12)
\path(87,237)(1737,237)
\path(912,1437)(912,87)
\path(12,87)(1512,837)
\path(312,837)(1812,87)
\path(912,1137)(762,912)
\path(803.603,1028.487)(762.000,912.000)(853.526,995.205)
\path(312,237)(612,237)
\path(492.000,207.000)(612.000,237.000)(492.000,267.000)
\path(1512,237)(1362,462)
\path(1453.526,378.795)(1362.000,462.000)(1403.603,345.513)
\path(912,237)(912,462)
\path(942.000,342.000)(912.000,462.000)(882.000,342.000)
\path(912,537)(686,654)
\path(806.359,625.472)(686.000,654.000)(778.774,572.189)
\path(1212,687)(1011,585)
\path(1104.434,666.056)(1011.000,585.000)(1131.586,612.551)
\path(613,461)(621,545)
\path(613,450)(561,526)
\end{picture}
}

\end{center}
\caption{The Fano plane, and a bijection of points and lines (indicated
by arrows from points to the corresponding line). \label{fano}}
\end{figure}

{\bf Counterexample 1.}  Eriksson and Linusson defined their
permutation array varieties over $\mathbb{C}$, so we begin with a
counterexample to realizability over $K = \mathbb{C}$, and it may be
read simply as an admonition to always consider a more general base
field (or indeed to work over the integers).  The Fano plane is the
projective plane over the field $\mathbb{F}_2$, consisting of $7$
lines $\ell_1$, \dots, $\ell_7$ and $7$ points $p_1$, \dots, $p_7$.
We may name them so that $p_i$ lies on $\ell_i$, as in
Figure~\ref{fano}.  Thus we have a configuration of $7$ flags over
$\mathbb{F}_2$.  (This is a projective picture, so this configuration
is in affine dimension $n=3$, and the points $p_i$ should be
interpreted as one-dimensional subspaces, and the lines $\ell_j$ as
two-dimensional subspaces, of $K^3$.)  The proof of
Theorem~\ref{t:EL-strata} is independent of the base field, so the
table of intersection dimensions of the flags yields a permutation
array.  However, a classical and straightforward argument in
projective coordinates shows that the configuration of
Figure~\ref{fano} may not be achieved over the complex numbers (or
indeed over any field of characteristic not $2$).  In particular, this
permutation array variety is not realizable over $\mathbb{C}$.  In
order to patch this counterexample, one might now restate the
Realizability Conjecture~\ref{c:realizability} by saying that there
always exists a field such that $X_{P}^{o}$ is nonempty.  However, the
problems have only just begun.

\begin{figure}[ht]
%\begin{center}
%\include{c3}
\hspace{-3.3in}\PSbox{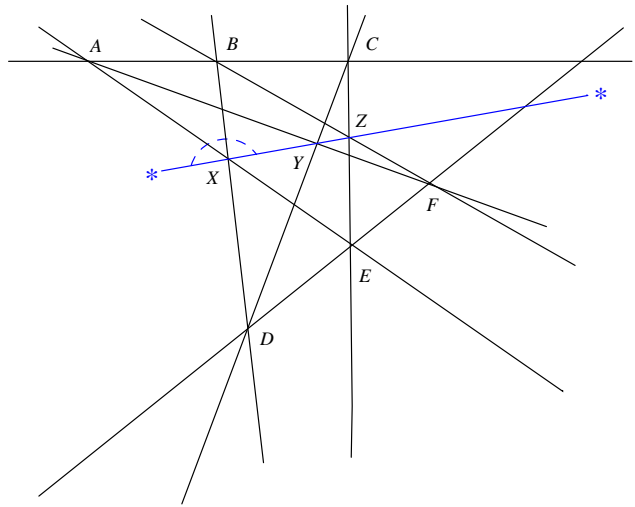}{0pt}{200pt}
%\end{center}
\caption{\label{counterexamplepappus}
Pappus' Theorem, and a counterexample to Realizability in dimension $d=3$ 
with $n=9$}
\end{figure}

{\bf Counterexample 2.}  We next sketch an elementary counterexample
for $n=3$ and $d=9$, over an arbitrary field, with the disadvantage
that it requires a computer check. Recall Pappus' Theorem in classical
geometry: if $A$, $B$, and $C$ are collinear, and $D$, $E$, and $F$
are collinear, and $X = AE \cap BD$, $Y= AF \cap CD$, and $Z= BF \cap
CE$, then $X$, $Y$, and $Z$ are collinear \cite[\S 3.5]{cg}.  The
result holds over any field.  A picture is shown in
Figure~\ref{counterexamplepappus}. (Ignore the dashed arc and the
stars for now.)

We construct an unrealizable permutation array as follows.  We imagine
that line $YZ$ does not meet $X$.  (In the figure, the
starred line $YZ$ ``hops over'' the point marked $X$.)  We construct a
counterexample with nine flags by letting the flags correspond to the
nine lines of our ``deformed Pappus configuration'', choosing points
on the lines arbitrarily.  We then construct the rank table of this
configuration, and verify that this corresponds to a valid permutation
array.  (This last step was done by computer.)  This
permutation array is not realizable, by Pappus' theorem.

%The Pappus configuration corresponds to a hyperplane intersecting 9
%flags in $\mathbb{R}^{3}$.  A permutation array corresponding to a
%near-Pappus arrangement, where the $*-*$ line ``hops'' over the real
%intersection is given by
%
%$$
%\begin{array}{cccc}
%(1 3 3 2 3 2 3 3 3 )&(2 3 3 3 1 3 2 3 3 )&(2 3 3 3 3 3 3 1 2) 
%
%&       (3 3 1 3 2 3 3 2 3) \\
%(3 3 2 3 3 1 3 3 2)&(3 3 2 2 3 3 1 3 3) 
%&       
%(3 1 3 3 2 2 3 3 3 )&(3 2 3 1 3 3 3 2 3)\\
%(3 3 3 3 3 3 2 3 1) 
%&       
%(2 3 2 3 3 3 3 3 3 )&(2 2 3 3 3 3 3 3 3)&(3 2 2 3 3 3 3 3 3) 
%\\
%(3 3 3 2 2 3 3 3 3)&(3 3 3 3 2 3 3 3 2)&(3 3 3 3 3 2 3 2 3 ) 
%&       
%(3 3 3 3 3 3 2 2 3)\\
%(3 2 3 3 3 3 2 3 3 )&(3 2 3 3 3 3 3 3 2) 
%&       
%(3 3 3 2 3 3 3 3 2)&(3 3 3 3 3 2 2 3 3)
%\end{array}
%$$

\begin{figure}[ht]
\begin{center}
\hspace{-2.1in}\PSbox{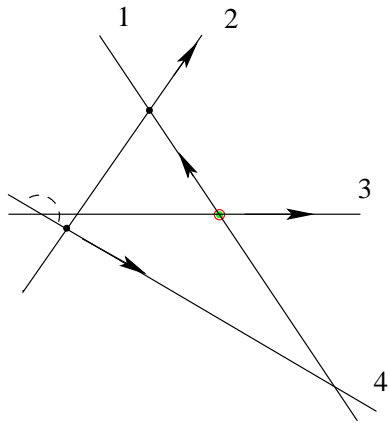}{0pt}{160pt} % was 60
\end{center}
\caption{\label{counterexampleplanar}
A counterexample to realizability with $n=d=4$}
\end{figure}

{\bf Counterexample 3.}  Our next example shows that realizability
already fails for $n=4$, $d=4$.  The projective intuition is as
follows.  Suppose $\ell_1$, $\ell_2$, $\ell_3$, $\ell_4$ are four
lines in projective space, no three meeting in a point, such that we
require $\ell_i$ and $\ell_j$ to meet, except (possibly) $\ell_3$ and
$\ell_4$.  This forces all 4 lines to be coplanar, so $\ell_3$ and
$\ell_4$ {\em must} meet.  Hence we construct an unrealizable
configuration as follows: we ``imagine'' (as in
Figure~\ref{counterexampleplanar}) that $\ell_3$ and $\ell_4$ don't meet.
Again, we must turn the projective picture in $\proj^3$ into linear
algebra in $4$-space, so the projective points in the figure correspond to
one-dimensional subspaces, the projective lines in the figure
correspond to two-dimensional subspaces of their respective flags, etc.
Again, the tail of each arrow corresponds with the point which lies on
the line the arrow follows. We construct the corresponding dot array:
$$\tableau{\emt & \emt & \emt & \emt \\
         \emt & \emt & \emt & \emt \\
         \emt & \emt & \emt & \emt \\
         \emt & 4 & \emt & \emt }
\hspace{.3in}
\tableau{\emt & \emt & \emt & \emt \\
         \emt & \emt & \emt & 4 \\
         \emt & \emt & \emt & \emt\\
         \emt & 4 & \emt & \cinum{3} }
\hspace{.3in}
\tableau{\emt & \emt & \emt & \emt \\
         \emt & \emt & \emt & 4 \\
         \emt & \emt & 4 & 3 \\
         \emt & 4 & 3 & 2 }
\hspace{.3in}
\tableau{\emt & \emt & \emt & 1 \\
         4 & \emt & \emt & \emt \\
         \emt & \emt & 3 & \emt \\
         \emt & 2 & \emt & \emt }
$$

\noindent Here the rows represent the flag $F^1_\ci$, columns
represent the flag $F^2_\ci$, numbers represent the flag $F^3_\ci$,
and the boards represent the flag $F^4_\ci$. This is readily checked
to be a permutation array.  The easiest way is to compare it to the
dot array for the ``legitimate'' configuration, where $F^3_2$ and
$F^4_2$ {\em do} meet, and using the fact that this second array is a
permutation array by Theorem~\ref{t:EL-strata}.  The only difference
between the permutation array above and the ``legitimate'' one is that
the circled 3 should be a 2.

%\vspace{1.5in}
%\hspace{.5in}\PSbox{planar.eps}{0pt}{60pt}

%
%The rank table for this collection of flags is represented
%by the permutation array:
%\bigskip
%$P=
%\tableau{\emt & \emt & \emt & \emt \\
%         \emt & \emt & \emt & \emt \\
%         \emt & \emt & \emt & \emt \\
%         \emt & 4 & \emt & \emt }
%\hspace{.3in}
%\tableau{\emt & \emt & \emt & \emt \\
%         \emt & \emt & \emt & \cinum{4} \\
%         \emt & \emt & \emt & \emt \\
%         \emt & 4 & \emt & 2 }
%\hspace{.3in}
%\tableau{\emt & \emt & \emt & \emt \\
%         \emt & \emt & \emt & 4 \\
%         \emt & \emt & 4 & 3 \\
%         \emt & 4 & 3 & 2 }
%\hspace{.3in}
%\tableau{\emt & \emt & \emt & 1 \\
%         4 & \emt & \emt & \emt \\
%         \emt & \emt & 3 & \emt \\
%         \emt & 2 & \emt & \emt }
%$

{\em Remark.}  Eriksson and Linusson have verified the Realizability
Conjecture~\ref{c:realizability} for $n=3$ and $d=4$ \cite[\S 3.1]{ELdecomposition}.  Hence
the only four open cases left are $n=3$ and $5 \leq d \leq 8$.  These
cases seem simple, as they involve (projectively) between $5$ and
$8$ lines in the plane.  Can these remaining cases be settled?

{\bf Further pathologies from Mn\"ev's universality theorem: failure
  of irreducible and equidimensionality.}  Mn\"ev's universality
theorem shows that permutation array schemes will be ``arbitrarily''
badly behaved in general, even for $n=3$.  Informally, Mn\"ev's
theorem states that given any singularity type of finite type over the
integers there is a configuration of projective lines in the plane
such that the corresponding {\em permutation array scheme} has that
singularity type.  By a singularity type of finite type over the
integers, we mean up to smooth parameters, any singularity cut out by
polynomials with integer co-efficients in a finite number of
variables.  See \cite{mnev1, mnev2} for the original sources, and
\cite[\S 3]{VakilC} for a precise statement and for an exposition of
the version we need.  (Mn\"ev's theorem is usually stated in a
different language of course.)

In particular, (i) permutation array schemes
need not be irreducible, answering a question raised 
in \cite[\S 1.2.3]{ELdecomposition}.  They can have arbitrarily
many components, indeed of arbitrarily many different dimensions.
(ii) Permutation array schemes need not be reduced, i.e.\ they
have genuine scheme-theoretic (or infinitesimal) structure not
present in the variety.  In other words, the definition
of permutation array schemes is indeed different from that
of permutation array varieties, and the equations \eqref{e:EL-strata}
do not cut out the permutation array varieties scheme-theoretically.
(iii) Permutation array schemes need not be equidimensional.
Hence the hope that permutation array varieties/schemes might
be well-behaved is misplaced.  In particular, the notion
of Bruhat order is problematic.  We suspect, for example, that
there exist two permutation array schemes $X$ and $Y$ such
that $Y$ is reducible, and some but not components of $Y$ lie
in the closure of $X$.

Although Mn\"ev's theorem is constructive, we have not attempted to
explicitly produce a reducible or non-reduced permutation array
scheme.

\section{Intersecting Schubert varieties}\label{s:intersections}

In this section, we consider a Schubert problem in $\flags$ of
the form
\[
X=X_{w^{1}}(E^{1}_{\ci})\cap X_{w^{2}}(E^{2}_{\ci}) \cap \dotsb
\cap X_{w^{d}}(E^{d}_{\ci})
\]
with $E^{1}_{\ci},\dots , E^{d}_{\ci}$ chosen generally and $\sum_i
\l(w^{i}) = \chs{n}{2}$.  We show there is a unique permutation array
$P$ for this problem if $X$ is nonempty, and we identify it.  In
Theorem~\ref{t:equations} we show how to use $P$ to write down
equations for $X$.  These equations can also be used to determine if
$E^{1}_{\ci},\dots , E^{d}_{\ci}$ are sufficiently general for
computing intersection numbers.
%For our purposes, those
%configurations of flags with the expected  number of solutions
%will be general enough.   [deleted because of previous sentence]
The number of solutions will always be either infinite or no greater than the
expected number.  The expected number is achieved on a dense open
subset of $\flags^d$.  It may be useful for the reader to refer to the
examples in Section~\ref{s:triple} while reading this section.

\begin{theorem}\label{t:unique}
If $X$ is 0-dimensional and nonempty, there exists a unique
permutation array $P \subset [n]^{d+1}$ such that 
%\begin{equation}\label{e:uniqueranktable}
$$  
\dim \left(E_{x_{1}}^{1}\cap E_{x_{2}}^{2}\cap
\cdots \cap E_{x_{d}}^{d} \cap F_{x_{d+1}}\right) = \mathrm{rk}P[x]
$$  
%\end{equation}
for all $F_\ci \in X$ and all $x\in [n]^{d+1}$.  Hence, $X$ is equal
to the fiber permutation array variety $X^{o}_{P}(E^{1}_{\ci},\dots,
E^{d}_{\ci})$.
\end{theorem}

As remarked in Section~\ref{s:realizability}, the projection of $P$
onto the first $d$ coordinates must be the transverse permutation
array $\transverse_{n,d}$.  We will prove the theorem by explicitly
constructing $P$.  As an immediate consequence, as the permutation
array corresponding to $d$ generally chosen flags $E^{1}_{\ci},\dots ,
E^{d}_{\ci}$ is given by $\transverse_{n,d}$, we have the following.

\begin{corollary}\label{c:empty}
If $P_{n} \neq \transverse_{n,d}$, then $X$ is the empty set.
\end{corollary}

When $d=4$, this corollary can often be used to detect when the
coefficients $c_{u,v}^{w}$ are zero in Equation~\ref{e:coefs}.  This
criterion catches $7$ of the $8$ zero coefficients in 3 dimensions,
$373$ of the $425$ in 4 dimensions, and $28920$ of the $33265$ in
dimension $5$.  The dimension $3$ case missed by this criterion is
presumably typical of what the criterion fails to see: there are no
$2$-planes in $3$-space containing three general $1$-dimensional
subspaces.  However, given a $2$-plane $V$, three general flags with
$1$-subspaces contained in $V$ are indeed transverse.

The corollary is efficient to apply.  For example, consider the
following three permutations\footnote{Anagrams of the name ``Richard
P. Stanley''.} in $S_{15}$:
\begin{itemize}
\item $u=$ A Children's Party
\item $v=$ Hip Trendy Rascal
\item $w =$ Raid Ranch Let Spy.
\end{itemize}
Using a computer, we can easily compute $P_{n}$ corresponding to
$X=X_u \cap X_v \cap X_{w}$:
\begin{tiny}
\[
P_{15}=\left[
\begin{array}{ccccccccccccccc}
  &&&&&&&&&&&&&& 11 \\
  &&&&&&&&&&&&& 14 &\\
  &&&&&&&&&&&& 14 & 11 & 3 \\
  &&&&&&&&&&& 15 &&&\\
  &&&&&&&&&& 15 & 14 & 11 & 5 &\\
  &&&&&&&&& 15 & 14 & 11 & 6 &&\\
  &&&&&&&& 15 & 8 &&&&&\\
  &&&&&&& 15 & 14 && 13 &&&&\\
  &&&&&& 15 & 14 & 8 & 4 &&&&&\\
  &&&&& 15 & 14 & 13 &&& 12 &&&&\\
  &&&& 15 & 14 & 13 & 8 & 3 &&&&&& 2 \\
  &&& 15 & 14 & 13 & 12 &&&& 11 & 9 &&&\\
  && 15 & 8 &&&& 2 &&&&&&& 1 \\
  & 15 & 13 &&& 11 &&&&& 10 &&&&\\
   15 & 13 & 8 & 7 &&&&&&&&&&&
\end{array} \right]
\]
\end{tiny}
Clearly, $P_{15} \neq T_{15,3}$ so $X$ is empty and $c_{u,v}^{w_{0}w}=0$.

%\begin{remark}\label{r:one.dim}
{\em Remark.}  The array $\transverse_{n,d}$ is an antichain under the dominance order on
$[n]^{d}$ so each element corresponds to a 1-dimensional vector space.
These lines will provide a ``skeleton'' for the given Schubert problem.
%\end{remark}

The proof of Theorem~\ref{t:unique} follows directly from the next 
lemma.

%\begin{proof}
%Let $V$ be an $n$-dimensional vector space.  For any two subspaces
%$U,W \subset V$ in general position, we have the formula
%$\dim\left(U \cap W\right)=\dim\left(U \right)
%+\dim\left(W \right) -n$.  So, assume  by induction 
%that $\dim \left(E_{x_{1}}^{1}\cap E_{x_{2}}^{2}\cap
%\cdots \cap E_{x_{d-1}}^{d-1}\right) = x_{1}+ \cdots +x_{d-1}-(d-2)n$, then 
%\begin{align*}\label{e:gnd.induction}
%\dim \left(E_{x_{1}}^{1}\cap E_{x_{2}}^{2}\cap
%\cdots \cap E_{x_{d}}^{d}\right) &= 
%\dim \left(E_{x_{1}}^{1}\cap E_{x_{2}}^{2}\cap
%\cdots \cap E_{x_{d-1}}^{d-1}\right) + 
%\dim \left(E_{x_{d}}^{d}\right) -n \\
%&= x_{1}+\dotsb +x_{d-1}-(d-2)n + x_{d}-n\\
%&= x_{1}+\dotsb +x_{d}-(d-1)n.
%\end{align*}
%Therefore, $\transverse_{n,d}$ is precisely the collection of indices
%$(x_{1},\dots ,x_{d})$ such that 
%\[
%\dim
%\left(E_{x_{1}}^{1}\cap E_{x_{2}}^{2}\cap \cdots \cap
%E_{x_{d}}^{d}\right) =1.
%\]
%By Theorem~\ref{t:EL-strata}, there exists a unique permutation array
%$P$ whose rank table agrees with the intersection data for
%$E^{1}_{\ci},\dots , E^{d}_{\ci}$.  We must have $\transverse_{n,d} \subset
%P$.  If there were an $x \in P \setminus \transverse_{n,d}$ then $x$ would be
%redundant, contradicting the definition of a permutation array.  Hence
%$P=\transverse_{n,d}$.
%\end{proof}

\begin{lemma}\label{l:recurrence}
Let $E^{1}_{\ci},\dots , E^{d}_{\ci}$ be generally chosen flags.
Let $w^{1},\dots , w^{d}$ be permutations in $S_{n}$ such that $\sum
\l(w^i) = \dim(\flags)$.  Let $F_{\ci}$ be a flag such that
$\mathrm{pos}(E^{i}_{\ci},F_{\ci})=w^{i}$ for each $1\leq i\leq d$.
Then the rank table of intersection dimensions among the components
of the $d+1$ flags is determined by the recurrence

$$
%\begin{equation}\label{e:recurrence}
\dim \left(E_{x_{s_{1}}}^{s_{1}}\cap E_{x_{s_{2}}}^{s_{2}}\cap 
\cdots \cap E_{x_{s_{k}}}^{s_{k}} \cap F_{j}\right)  \hspace{2in}
%\end{equation}
$$
$$
%\begin{equation*}
\hspace{1in}= 
\mathrm{max}
\begin{cases}
 \dim 
\left(E_{x_{s_{1}}}^{s_{1}}\cap F_{j}\right) + 
\dim \left( E_{x_{s_{2}}}^{s_{2}}\cap
\cdots \cap E_{x_{s_{k}}}^{s_{k}} \cap F_{j}\right)   -j\\
 \dim\left(E_{x_{s_{1}}}^{s_{1}}\cap E_{x_{s_{2}}}^{s_{2}}\cap
\cdots \cap E_{x_{s_{k}}}^{s_{k}} \cap F_{j-1}\right)
\end{cases}
%\end{equation*}
$$
where $1\leq s_{1}<s_{2}<\dotsb <s_{k}\leq d$,\, 
$k\geq 2$, \,
 $1 \leq x_{s_{i}}
\leq n-1$ for each $1\leq i\leq k$,
$
\dim\left(F_{0}\right) =0$, and $\dim
\left(E_{i}^{a}\cap F_{j}\right)$ is determined by the rank table
corresponding to the permutation $w^{a}$.
\end{lemma}

Note that this recurrence determines the full intersection table since
\[
\dim \left(E_{x_{1}}^{1}\cap E_{x_{2}}^{2}\cap
\cdots \cap E_{x_{d}}^{d} \cap F_{j}\right)
= \dim \left(E_{x_{s_{1}}}^{s_{1}}\cap E_{x_{s_{2}}}^{s_{2}}\cap 
\cdots \cap E_{x_{s_{k}}}^{s_{k}} \cap F_{j}\right)
\]
if $x_{i}=n$ for each $i \in [d] \setminus \{s_{1},\dots , s_{k} \}$. 

\begin{proof}
Set $U=E_{x_{s_{1}}}^{s_{1}} \cap F_{j}$ and
$V=E_{x_{s_{2}}}^{s_{2}}\cap \cdots \cap E_{x_{s_{k}}}^{s_{k}} \cap
F_{j}$.  Since $F_{j-1} \subset F_{j}$ and $\textrm{dim}(F_{j})=
1+\textrm{dim}(F_{j-1})$, we have
\[
\dim (U \cap V \cap F_{j-1}) \leq \dim (U \cap V )
\leq 1+ \dim (U \cap V \cap F_{j-1}).
\]
We also know that $\dim (U \cap V )\geq \dim (U ) +
\dim (V) - j$ since $U,V \subset F_{j}$ and
$\dim(F_{j})=j$.  We need to show $\dim (U \cap V \cap
F_{j-1}) < \dim (U \cap V )$ if and only if $\dim (U
\cap V) = \dim (U ) + \dim (V) - j$.  

Let $W=U \cap V \cap F_{j-1}$ and choose $U'$ and $V'$ so that $U=W\oplus U'$ and $V=W\oplus
V'$.  By the assumption that $E^{1}_{\ci},\dots , E^{d}_{\ci}$
are general we have 
%todo{Ravi, our new definition of general doesn't exactly cover this. 
%What should we do? -- S.  I've cut our new definition of general! --- R}
\[
\dim(U'\cap V') = \mathrm{max}
\begin{cases}
\dim(U')+ \dim(V') - \left(j- \dim(W)
\right)\\
0  
\end{cases}.
\]
Therefore, $\dim (U \cap V )>\dim(W)$ if and only if 
$\dim(U'\cap V')>0 $ if and only if  $\left(\dim(U')+
\dim(V') - \left(j- \dim(W)\right) \right) >0 $ if and only if 
\begin{align*}
\dim(U\cap V) &= \dim(W) + \dim(U' \cap V')\\
&= \dim(W) + \dim(U')+
\dim(V') - \left(j- \dim(W)
\right) \\
&=\dim (U ) + \dim (V) - j.
\end{align*}

\end{proof}

\begin{theorem}\label{t:equations}
Let $ X=X_{w^{1}}(E^{1}_{\ci})\cap X_{w^{2}}(E^{2}_{\ci}) \cap \dotsb
\cap X_{w^{d}}(E^{d}_{\ci}) $ be a 0-dimensional intersection, with
$E^{1}_{\ci},\dots , E^{d}_{\ci}$ general.  Let $P \subset
[n]^{d+1}$ be the unique permutation array associated to this
intersection by the recurrence in Lemma~\ref{l:recurrence}.  Let
$V(E^{1}_{\ci},\dots , E^{d}_{\ci})= \{v_{x} \given x \in
\transverse_{n,d} \}$ be a collection of vectors chosen such that
$v_{x} \in E_{x_{1}}^{1}\cap E_{x_{2}}^{2}\cap \cdots \cap
E_{x_{d}}^{d}$.  Then polynomial equations defining $X$ can be
determined simply by knowing $P$ and $V(E^{1}_{\ci},\dots ,
E^{d}_{\ci})$.
\end{theorem}

\begin{proof}
Given $P\in [n]^{d+1}$, let $P_{1},\dots, P_{n}$ be the sequence of
permutation arrays in $[n]^{d}$ defined by the EL-algorithm in
Theorem~\ref{t:EL-algorithm}.  If $F_{\ci} \in X$, then by
Corollary~\ref{c:checkerboards} $P_{i}$ is the unique permutation
array encoding $\dim(E_{x_{1}}^{1}\cap E_{x_{2}}^{2}\cap \cdots \cap
E_{x_{d}}^{d} \cap F_{i})$.  Furthermore, for each $x \in P_{i}, 1\leq
i\leq n$, we could choose a representative vector in the corresponding
intersection, say $v_{x}^{i} \in E_{x_{1}}^{1}\cap
E_{x_{2}}^{2}\cap \cdots \cap E_{x_{d}}^{d} \cap F_{i}$.  Define
\begin{align*}
&V_{i}= \{v_{y}^{i} \given y \in P_{i} \} \\
&V_{i}[x] = \{v_{y}^{i} \given y \in P_{i}[x] \}.
\end{align*}
In fact, we can choose the vectors $v_{x}^{i}$ so that $v_{x}^{i} \not
\in \mathrm{Span}(V_{i}[x] \setminus \{v_{x}^{i} \})$ since the rank
function must increase at position $x$.  Therefore, we would have
\begin{equation}\label{e:rank.equations.1} % used
\mathrm{v.rk}(V_{i}[x]) = \mathrm{rk}(P_{i}[x])
\end{equation}
for all $x \in [n]^{d}$ and all $1\leq i\leq n$ where
$\mathrm{v.rk}(S)$ is the dimension of the vector space spanned by
the vectors in $S$.
These rank equations define $X$.  

Let $V_{n}=V(E^{1}_{\ci},\dots ,
E^{d}_{\ci})$ be the finite collection of vectors in the case $i=n$.
Given $V_{i+1}$, $P_{i+1}$ and $P_{i}$, we compute
\[
V_{i}
= \{v_{x}^{i} \given x \in P_{i} \}
\]
recursively as follows.  If $x \in P_{i} \cap P_{i+1}$ then set
$$
v_{x}^{i}=v_{x}^{i+1}.  
$$
If $x \in P_{i} \setminus P_{i+1}$ and $y,\dots , z$ is a \textit{basis
set} for $P_{i+1}[x]$, i.e.\ $v_{y}^{i+1},\dots, v_{z}^{i+1}$ are
independent and span the vector space generated by all $v_{w}^{i+1}$
with $w \in P_{i+1}[x]$, then set
\begin{equation}\label{e:vectors} % used
v_{x}^{i} = c^{i}_{y}v_{y}^{i}+\dotsb + c^{i}_{z}v_{z}^{i}
\end{equation}
where $c^{i}_{y}, \dots, c^{i}_{z}$ are indeterminate.  Now the same
rank equations as in \eqref{e:rank.equations.1} must hold.  In fact, it
is sufficient in a 0-dimensional variety $X$ to require only
\begin{equation}\label{e:rank.equations.2} % used
\mathrm{v.rk}\{v_{y}^{i} \given y \in P_{i}[x] \} \leq  \mathrm{rk}(P_{i}[x])
\end{equation}
for all $x \in [n]^{d}$ and all $1\leq i\leq n$.  Let
$\mathrm{minors}_{k}(M)$ be the set of all $k \times k$ determinantal
minors of a matrix M.  Let $M(V_{i}[x])$ be the matrix whose rows are
given by the vectors in $V_{i}[x]$.

Then, the equations \eqref{e:rank.equations.2}
can be rephrased as
\begin{equation}\label{e:rank.equations.3} % used
\mathrm{minors}_{rk(P_{i}[x])+1}(M(V_{i}[x]))=0  
\end{equation}
for all $1\leq i<n$ and $x \in [n]^{d}$ such that $\sum
x_{i}>(d-1)n$.

For each set of solutions $S$ to the equations in
\eqref{e:rank.equations.3}, we obtain a collection of vector sets by
substituting solutions for the indeterminates in the formulas
\eqref{e:vectors} for the vectors.  We can further eliminate variables
whenever a vector depends only on one variable $c_{y}^{i}$ by setting
it equal to any nonzero value which does not force another
$c_{z}^{j}=0$.  If ever a solution implies $c_{z}^{j}=0$, then the
choice of $E^{1}_{\ci},\dots , E^{d}_{\ci}$ was not general.  Let
$V_{1}^{S},\dots V_{n}^{S}$ be the final collection of vector sets
depending on the solutions $S$.  Since $X$ is 0-dimensional, if
$V_{1}^{S},\dots V_{n}^{S}$ depends on any indeterminate then
$E^{1}_{\ci},\dots , E^{d}_{\ci}$ was not general.  Let $F^{S}_{i}$ be
the span of the vectors in $V_{i}^{S}$.  Then the flag
$F_{\ci}^{S}=(F_{1}^{S},\dots, F_{n}^{S})$ satisfies all the rank
conditions defining $X=X_{P}^{o}(E^{1}_{\ci},\dots, E^{d}_{\ci})$.
Hence, $F_{\ci}^{S} \in X$.
\end{proof}

%\begin{remark}\label{r:improvements}
{\em Remark.}  There are too many equations and indeterminates
involved in the equations \eqref{e:rank.equations.3} to solve this
system simultaneously in practice.  First, it is useful to solve all
equations pertaining to $V_{i+1}$ before computing the initial form of
the vectors in $V_{i}$.  Second, we have found that proceeding through
all $x \in [n]^{d}$ such that $\sum x_{i}>(d-1)n$ in lexicographic
order works well, with the additional caveat that if $P_{i}[x]=\{x \}$
then the rank of the matrix $M$ with rows determined by $\{x \}\cup
(P_{i}\cap P_{i+1})$ must have rank at most $i$.  Solve all of the
determinantal equations implying the rank condition
$\mathrm{v.rk}(V_{i}[x]) = \mathrm{rk}(P_{i}[x])$ simultaneously and
substitute each solution back into the collection of vectors before
considering the next rank condition.

%\end{remark}

\begin{corollary}\label{c:generic}
The equations appearing in \eqref{e:rank.equations.3} provide a test
for determining if $E^{1}_{\ci},\dots , E^{d}_{\ci}$ is sufficiently
general for the given Schubert problem.  Namely, the number of flags
satisfying the equations \eqref{e:rank.equations.3} is the generic
intersection number if each indeterminate $c_{z}^{i} \neq 0$ and the
solution space determined by the equations is 0-dimensional.
\end{corollary}

%todo{This corollary depends on our definition of generic.  Here,
%generic must mean that the number of solutions equals the generic
%number of solutions.  It does not imply that there are no hidden
%dependencies.  -S}

\section{The key example: Triple intersections}\label{s:triple}

We now implement the algorithm of the previous section in an important special
case.
Our goal is to describe a method for directly
identifying all flags in $X=X_{u}(E^{1}_{\ci}) \cap X_{v}(E^{2}_{\ci}) \cap
X_{w}(E^{3}_{\ci}) $ when $\l (u) + \l (v) + \l (w) = \chs{n}{2}$ and
$E^{1}_{\ci}$, $E^{2}_{\ci}$, and $E^{3}_{\ci}$ are in general position.  This
gives a method for computing the structure constants in the cohomology
ring of the flag variety from equations \eqref{e:coefs} and
\eqref{e:triple-intersection} .

There are two parts to this algorithm.  First, we use the recurrence
of  Lemma~\ref{l:recurrence} to find the unique
permutation array $P\subset [n]^{4}$ with position vector $(u,v,w)$
such that $P_{n}=\transverse_{n,3}$.  Second, given $P$ we use the equations in
\eqref{e:rank.equations.3} to find all flags in $X$.  

As a demonstration, we explicitly compute the flags in $X$ in two
cases.  For convenience, we work over $\mathbb{C}$, but of course the
algorithm is independent of the field.  In the first there is just one
solution which is relatively easy to see ``by eye''.  In the second
case, there are two solutions, and the equations are more complicated.
The algorithm has been implemented in Maple and works well on examples
where $n\leq 8$.  

\begin{example}\label{ex:flags.1}  \end{example}
Let $u=(1,3,2,4)$, $v=(3,2,1,4)$, $w=(1,3,4,2)$.  The sum of their
lengths is $1+3+2=6 = \chs{n}{2}$.  The unique permutation array $P
\in [4]^{4}$ determined by the recurrence in Lemma~\ref{l:recurrence}
consists of the following dots:
\[
\begin{array}{ccccc}
(4 4 2 1)&       (4 1 4 2) &    (2 4 4 2)&       (4 2 3 3) &    (3 2 4 3)\\
 (3 4 3 3)&      (4 4 1 4)&      (4 3 2 4)&      (3 4 2 4) &    (3 3 3 4)\\
 (2 4 3 4) &    (2 3 4 4) &     (1 4 4 4)
\end{array}
\]
The EL-algorithm produces the following list of permutation arrays
$P_{1},P_{2},P_{3},P_{4}$ in $[4]^{3}$ corresponding to $P$:

\[
\begin{array}{cccc}
\tableau{\emt & \emt & \emt & \emt \\
         \emt & \emt & \emt & \emt \\
         \emt & \emt & \emt & \emt \\
         \emt & \emt & \emt & 2 }
&
\tableau{\emt & \emt & \emt & \emt \\
         \emt & \emt & \emt & 4 \\
         \emt & \emt & \emt & \emt \\
         4 & \emt & \emt & 2 }
&
\tableau{\emt & \emt & \emt & \emt \\
         \emt & \emt & \emt & 4 \\
         \emt & 4        & \emt & 3 \\
         4        & 3        & \emt & 2 }
&
\tableau{\emt & \emt & \emt & 4 \\
         \emt & \emt & 4 & 3 \\
         \emt & 4 & 3 & 2 \\
         4 & 3 & 2 & 1 }
\\
P_{1} & P_{2} & P_{3} & P_{4} 
\end{array}
\]

\noindent 
Notice that $P_{4}$ is the transverse permutation array
$\transverse_{4,3}$. Notice also how to read $u$, $v$, and $w$ from
$P_1$, \dots, $P_4$: $P_i$ has one less row than $P_{i+1}$; listing
these excised rows from right to left yields $u$.  Similarly, listing
the excised {\em columns} from right to left yields $v$, and listing
the excised {\em numbers} from right to left yields $w$ (see the
example immediately above).

% Also, the permutation $u$ determines the
%order in which the rows are eliminated reading right to left, $v$
%determines the order for the columns, and $w$ determines the order for
%the numbers.

We want to specify three transverse fixed flags $E^{1}_{\ci}$, $E^{2}_{\ci}$,
$E^{3}_{\ci}$.  It will be notationally
convenient to represent a vector $v=(v_{1}, \dots , v_{n})$ by the
polynomial $v_{1} + v_{2} x + \cdots + v_{n} x^{n-1}$.  We choose three
flags, or equivalently three ``transverse'' ordered bases,
as follows:
$$
%\begin{equation}\label{e:transverse.flags}
\begin{array}{ccl}
E^{1}_{\ci} &=& \langle 1,x,x^{2},x^{3} \rangle \\
E^{2}_{\ci} &=& \langle x^{3}, x^{2},x ,1 \rangle \\
E^{3}_{\ci} &=& \langle (x+1)^{3}, (x+1)^{2}, (x+1), 1 \rangle 
\end{array}
%\end{equation}
$$

We will show that the only flag in $X_{u}(E^{1}_{\ci})\cap
X_{v}(E^{2}_{\ci})\cap X_{w}(E^{3}_{\ci})$ is 
\begin{equation}\label{e:eq1}
F_{\ci}=\langle 2+3x-x^{3} , \, x^{3}, \, x^{2}, \, 1
\rangle .
\end{equation}

For each element $(i,j,k)$ in $P_{4}$, we choose a vector in the
corresponding 1-dimensional intersection $E^{1}_{i}\cap E^{2}_{j}\cap
E^{3}_{k}\cap F_{4}$ and put it in position $(i,j)$ in the matrix
below:
$$V(E^{1}_{\ci}, E^2_{\ci}, E^{3}_{\ci}) = V_{4} = 
\left[
\begin{array}{ccccc}
0&      0&      0&      1\\
0&      0&      x&      x+1\\
0&      x^{2}&  x(x+1)& (x+1)^{2}\\
x^{3}&  x^{2}(x+1)&x(x+1)^{2}&  (x+1)^{3}
\end{array} \right].
$$
In $P_{3}$, every element in the 4th column is covered by a subset in
the antichain removed from $P_{4}$.  This column adds only one degree
of freedom so we establish $V_{3}$ by adding only one variable in
position $(2,4)$ and solving all other rank two equations in terms of
this one:
$$V_{3} = 
\left[
\begin{array}{cccc}
0&      0&      0&     0\\
0&      0&      0&      (1+ x)  + c x\\
0&      x^{2}&  0&      1+x +(1+c)x(1+x)\\
x^{3}&  x^{2}(x+1)&0&  (x+1)^{2} + c x(x+1)^{2}
\end{array}
\right].
$$
According to equation~\eqref{e:vectors} the entry in position
$(4,2)$ can have two indeterminates: $b(1+ x) + c x$, where $b, c \neq 0$.
As any two linearly dependent ordered pairs $(b,c)$ yield the same configuration
of subspaces, we may normalize $b$ to $1$.

Once $V_{3}$ is determined, we  find the vectors in
$V_{2}$. In $P_{2}$, every element is contained in $P_{3}$, so $V_{2}$
is a subset of $V_{3}$:
$$V_{2} = 
\left[
\begin{array}{cccc}
0&      0&      0&     0\\
0&      0&      0&      1+ (2+c) x\\
0&      0&      0&      0\\
x^{3}&  0&0&  (x+1)^{2} + c x(x+1)^{2}
\end{array}\right].
$$
The rank of $P_{2}$ is $2$, so all $3\times 3$ minors of the following
matrix must be zero:
\[
\left(
\begin{array}{cccc}
0 &     0 &     0 &     1\\
1 & 2+c & 0 & 0\\
1 &     2+c &    1+2c & c
\end{array}
  \right).
\]
In particular, $1+2c=0$, so the only solution is $c=-\frac{1}{2}$.
Substituting for $c$, we have
$$
V_{2}^{S} = 
\left[
\begin{array}{cccc}
0&      0&      0&     0\\
0&      0&      0&      1+ \frac{3}{2} x\\
0&      0&      0&      0\\
x^{3}&  0&0&  (x+1)^{2} - \frac{1}{2} x(x+1)^{2}
\end{array} \right].
$$
Finally $P_{1}$ is contained in $P_{2}$, so $V_{1}^{S}$ contains
just the vector 
\[
v^{1}_{(4,4,2)}=v^{2}_{(4,4,2)} = (x+1)^{2} - \frac{1}{2} x(x+1)^{2} =
\frac{1}{2}(2+3x-x^{3}).
\]

Therefore, there is just one solution, namely the flag spanned by the
collections of vectors $V_{1}^{S},V_{2}^{S},V_{3}^{S},V_{4}^{S}$ which
is equivalent to the flag in \eqref{e:eq1}.

If we choose an arbitrary general collection of three flags, we can
always change bases so that we have the following situation:
$$
%\begin{equation}\label{e:transverse.flags}
\begin{array}{ccl}
E^{1}_{\ci} &=& \langle 1,x,x^{2},x^{3} \rangle \\
E^{2}_{\ci} &=& \langle x^{3}, x^{2},x ,1 \rangle \\
E^{3}_{\ci} &=& \langle a_{1}+a_{2}x+a_{3}x^{2}+x^{3},
 \, b_{1}+b_{2}x+x^{2},\, c_{1}+x,\, 1\rangle
\end{array}
%\end{equation}
$$
Using these coordinates, the same procedure as above will produce the
unique solution 
$$F_{\ci}=\langle
(a_{1}-a_{3}b_{1})+(a_{2}-b_{2}a_{3})x - x^{3}, \, x^{3},\,
x^{2},\, 1
\rangle .$$

\begin{example}\label{ex:2} \end{example}
This example is of a Schubert problem with multiple solutions.
Let $u=(1, 3, 2, 5, 4, 6)$, $v=(3,5,1,2,4,6)$, $w=(3,1,6,5,4,2)$.  If
$P$ is the  unique permutation array in $[6]^{4}$ determined by the
recurrence in Lemma~\ref{l:recurrence} for $u,v,w$ then the
EL-algorithm produces the following list of permutation arrays
$P_{1},\dots,P_{6}$ in $[6]^{3}$ corresponding to $P$:
\[
\tableau{\emt& \emt& \emt& \emt& \emt& \emt\\
   \emt& \emt& \emt& \emt& \emt& \emt\\
   \emt& \emt& \emt& \emt& \emt& \emt\\
   \emt& \emt& \emt& \emt& \emt&\emt\\
   \emt& \emt& \emt& \emt& \emt& \emt\\
   \emt& \emt& \emt&\emt& \emt&2}
\hspace{.3in}
\tableau{\emt& \emt& \emt& \emt& \emt& \emt\\
   \emt& \emt& \emt& \emt& \emt& \emt\\
   \emt& \emt& \emt& \emt& \emt& \emt\\
   \emt& \emt& \emt& \emt& \emt&4\\
   \emt& \emt& \emt& \emt& \emt& \emt\\
   \emt& \emt& \emt&4& \emt&2}
\hspace{.3in}
\tableau{\emt& \emt& \emt& \emt& \emt& \emt\\
   \emt& \emt& \emt& \emt& \emt& \emt\\
   \emt& \emt& \emt& \emt& \emt& \emt\\
   \emt& \emt& \emt& \emt& \emt&4\\
   \emt& \emt& \emt&5& \emt& \emt\\
   \emt&5& \emt&4& \emt&2}
\hspace{.5in}
\]

\[\hspace{.5in}
\tableau{  
   \emt& \emt& \emt& \emt& \emt& \emt\\
   \emt& \emt& \emt& \emt& \emt&6\\
   \emt& \emt& \emt& \emt& \emt& \emt\\
   \emt& \emt& \emt&6& \emt&4\\
   \emt&6& \emt&5& \emt& \emt\\
  6&5& \emt&4& \emt&2&}
\hspace{.2in}
\tableau{  
   \emt& \emt& \emt& \emt& \emt& \emt\\
   \emt& \emt& \emt& \emt& \emt&6\\
   \emt& \emt& \emt& \emt&6&5\\
   \emt& \emt& \emt&6&5&4\\
   \emt&6& \emt&5&4&2\\
  6&5& \emt&4&2&1&}
\hspace{.2in}
\tableau{  
   \emt& \emt& \emt& \emt& \emt&6\\
   \emt& \emt& \emt& \emt&6&5\\
   \emt& \emt& \emt&6&5&4\\
   \emt& \emt&6&5&4&3\\
   \emt&6&5&4&3&2\\
  6&5&4&3&2&1&}
\]
We take the following triple of fixed flags:

\[
\begin{array}{ccl}
E^{1}_{\ci} &=& \langle 1, x, \dots , x^{5} \rangle \\
E^{2}_{\ci} &=& \langle x^{5}, \dots ,x,1 \rangle \\
E^{3}_{\ci} &=& \langle (1+x)^{5}, (1+x)^{4}, \dots , 1 \rangle 
\end{array}
\]
The third flag is clearly not  chosen generally but leads to two
solutions to this Schubert problem which is the generic number of
solutions.  We prefer to work with explicit but simple numbers here to
demonstrate the computation without making the formulas too
complicated.

The vector table associated to $P_{6}$ is easily determined by
Pascal's formula:
\begin{tiny}
\[
 \left[ \begin {array}{cccccc} []&[]&[]&[]&[]&[1,0,0,0,0,0]
\\\noalign{\medskip}[]&[]&[]&[]&[0,1,0,0,0,0]&[1,1,0,0,0,0]
\\\noalign{\medskip}[]&[]&[]&[0,0,1,0,0,0]&[0,1,1,0,0,0]&[1,2,1,0,0,0]
\\\noalign{\medskip}[]&[]&[0,0,0,1,0,0]&[0,0,1,1,0,0]&[0,1,2,1,0,0]&[1,
3,3,1,0,0]\\\noalign{\medskip}[]&[0,0,0,0,1,0]&[0,0,0,1,1,0]&[0,0,1,2,1
,0]&[0,1,3,3,1,0]&[1,4,6,4,1,0]\\\noalign{\medskip}[0,0,0,0,0,1]&[0,0,0
,0,1,1]&[0,0,0,1,2,1]&[0,0,1,3,3,1]&[0,1,4,6,4,1]&[1,5,10,10,5,1] \end
{array} \right]
\]
\end{tiny}
%%%%%%%%%%%%%%%%%%%%%%%%%%%%%%%%%%%%%%%%%%%%%%%%%%%%%%%%%%%%%%%%%%5

The vector table associated to $P_{5}$ has one degree of freedom.  The
vector in position $(3,5)$ is freely chosen to be $x+c\, x^{2}$.  Then
for all other points in $P_{5} \setminus P_{6}$ we can solve a rank 2
equation which determines the corresponding vector in terms of $c$.  Therefore, $V_{5}$ becomes 

\begin{tiny}
\[
 \left[ \begin {array}{cccccc}
[]&[]&[]&[]&[]&[]\\
\noalign{\medskip}[]& []&[]&[]&[]&[1,{\frac{10(c-1)}{3c}} ,0,0,0,0]\\
\noalign{\medskip}[]&[]&[]&[]&[0,1,c,0,0,0]&[1,{\frac {-13+10\,c}{3(c-1)}},{\frac {-10+7\,c}{3(c-1)}},0,0,0]\\
\noalign{\medskip}[]&[]&[]&[0,0,1,\frac{-6}{c-4 },0,0]&[0,1,{\frac {-4+7\,c}{2+c}},{\frac {6(c-1)}{2+c}},0,0]&[1,6\,{\frac {-5+3\,c}{-8+5\,c}}, {\frac {3(-12+7\,c)}{-8+5\,c}},{\frac {2(-7+4\,c)}{-8+5\,c}},0,0]\\
\noalign{\medskip}[]&x^4&[]&[0,0,1,\frac{-6}{c-4 },{\frac {-2-c}{c-4}},0]&[0,1,{\frac {3(-4+3\,c)}{2(c-1)}},{\frac {3(-3+2\,c)}{c-1}},{ \frac{-8+5\,c}{2(c-1)}},0]&[1,4,6,4,1,0]\\
\noalign{\medskip}x^{5}&x^4+x^{5}&[]&[0,0,1,\frac{-6}{c-4 },{\frac {-3c}{c-4}},{\frac{2(1-c)}{c-4}}]&[0,1,4,6,4,1]&[1,5,10,10,5,1]\end {array} \right]\]
\end{tiny}
%%%%%%%%%%%%%%%%%%%%%%%%%%%%%%%%%%%%%%%%%%%%%%%%%%%%%%%%%%%%%%%%%%

Every vector in $V_{4}$ appears in $V_{5}$, but now some of them are
subject to new rank conditions:
\begin{tiny}
\[
 \left[ \begin {array}{cccccc}
[]&[]&[]&[]&[]&[]\\
\noalign{\medskip}[]& []&[]&[]&[]&[1,{\frac{10(c-1)}{3c}} ,0,0,0,0]\\
\noalign{\medskip}[]&[]&[]&[]&[]&[]\\
\noalign{\medskip}[]&[]&[]&[0,0,1,\frac{-6}{c-4 },0,0]&[]&[1,{\frac {6(-5+3\,c)}{-8+5\,c}}, {\frac {3(-12+7\,c)}{-8+5\,c}},{\frac {2(-7+4\,c)}{-8+5\,c}},0,0]\\
\noalign{\medskip}[]&x^4&[]&[0,0,1,\frac{-6}{c-4 },-{\frac {2+c}{c-4}},0]&[]&[]\\
\noalign{\medskip}x^{5}&x^4+x^{5}&[]&[0,0,1,\frac{-6}{c-4 },{\frac {-3c}{c-4}},{\frac{-2c+2}{c-4}}]&[]&[1,4+d,6+4d,4+6d,1+4d,d]\end {array} \right].\]
\end{tiny}
%%%%%%%%%%%%%%%%%%%%%%%%%%%%%%%%%%%%%%%%%%%%%%%%%%%%%%%%%%%%%%%%%%
In particular, the top $3$ vectors should span a two-dimensional subspace.
This happens if the following matrix has rank $2$:
$$
\left[\begin {array}{cccccc}
1&{\frac {10(c-1)}{3c}}&0&0&0&0\\
1&{\frac {6(-5+3\,c)}{-8+5\,c}}&{\frac {3(-12+7\,c)}{-8+5\,c}}&{\frac {2(-7+4\,c)}{-8+5\,c}}&0&0
\\
0&0&1&\frac{-6}{c-4}&0&0
\end {array} \right] 
$$
or equivalently if the following nontrivial minors of the matrix are zero
\[
\left[{\frac {4(10\,c+{c}^{2}-20)}{3c \left( -8+5\,c \right) }},{\frac 
{-8(10\,c+{c}^{2}-20)}{ \left( -8+5\,c \right)  \left( c-4 \right) c}},
{\frac {-8(10\,c+{c}^{2}-20)}{ \left( -8+5\,c \right)  \left( c-4 \right) }
},{\frac { -8\left( c-1 \right)  \left( 10\,c+{c}^{2}-
20 \right) }{3 \left( -8+5\,c \right)  \left( c-4 \right) c}} \right].
\]
All rank $3$ minors will be zero if ${c}^{2}+ 10\, c -20 =0$, or $c=-5\pm
3 \, \sqrt{5}$.  Plugging each solution for $c$ into the vectors gives
the two solutions $V_{4}^{S_{1}}$ and $V_{4}^{S_{2}}$.  For example,
using $c=-5+ 3 \, \sqrt{5}$ and solving a single rank 2 equation involving $d$ gives:

%%%%%%%%%%%%%%%%%%%%%%%%%%%%%%%%%%%%%%%%%%%%%%%%%%%%%%%%%%%%%%%%%%
\begin{tiny}
\[ \left[ \begin {array}{cccccc} []&[]&[]&[]&[]&[]\\
\noalign{\medskip}[]&[]&[]&[]&[]&[1,10\,{\frac {2+\sqrt {5}}{5+3\,\sqrt {5}}},0,0,0,0]\\
\noalign{\medskip}[]&[]&[]&[]&[]&[]\\
\noalign{\medskip}[]&[]&[]&[0,0,1,\frac{2}{\left( 3+\sqrt {5} \right)},0,0]&[]&[1,2\,{\frac {20+9\,\sqrt {5}}{11+5\,\sqrt {5}}},{\frac {47+21\,\sqrt {5}}{11+5\,\sqrt {5}}},2\,{\frac {9+4\,\sqrt {5}}{11+5\,\sqrt {5}}},0,0]\\
\noalign{\medskip}[]&x^4&[]&[0,0,1,\frac{2}{\left( 3+\sqrt {5} \right)},-{\frac {1+\sqrt {5}}{3+\sqrt {5}}},0]&[]&[]\\
\noalign{\medskip}x^{5}&x^4+x^{5}&[]&[0,0,1,\frac{2}{\left( 3+\sqrt {5}
\right)},-{\frac {5+3\,\sqrt {5}}{3+\sqrt {5}}},{\frac
{-2(2+\sqrt {5})}{3+\sqrt {5}}}]&[]&[1,\frac{5+\sqrt {5}}{2},2\,\sqrt
{5},-5+3\,\sqrt {5},-5+2\,\sqrt {5},\frac{-3+\sqrt {5}}{2}]\end {array}
\right] \]
\end{tiny}
%%%%%%%%%%%%%%%%%%%%%%%%%%%%%%%%%%%%%%%%%%%%%%%%%%%%%%%%%%%%%%%%%%
The remaining vectors in $V_{1}^{S_{1}}, V_{2}^{S_{1}}, V_{3}^{S_{1}}$
will be a subset of $V_{4}^{S_{1}}$ so no further equations need to be
solved, and similarly for $V_{4}^{S_{2}}$.

%%%%%%%%%%%%%%%%%%%%%%%%%%%%%%%%%%%%%%%%%%%%%%%%%%%%%%%%%%%%%%%%%%

%\begin{verbatim}
%CL-USER(22): (print-schubert-triple-arrays '(1 3 2 5 4 6) '(3 5 1 2 4 6) 
%(mult-w0-left '(3 1 6 5 4 2) 'a))

%((5 5 1 0) (5 3 3 1) (3 5 3 1) (5 1 4 2) (4 3 4 2) (5 0 5 3) (4 1 5 3)
% (3 3 5 3) (1 5 5 3) (2 5 4 4) (5 4 1 4) (4 4 3 4) (3 4 4 4) (2 4 5 4)
% (5 5 0 4) (4 5 1 4) (0 5 5 5) (5 2 3 5) (4 2 4 5) (3 2 5 5) (5 3 2 5)
% (4 4 2 5) (3 5 2 5) (1 4 5 5) (2 3 5 5) (1 5 4 5) (2 4 4 5) (3 3 4 5)
% (2 5 3 5) (3 4 3 5) (4 3 3 5)) 
%    ......    ......    ......    ......    ......    .....5
%    ......    ......    ......    .....5    .....5    ....54
%    ......    ......    ......    ......    ....54    ...543
%    ......    .....3    .....3    ...5.3    ...543    ..5432
%    ......    ......    ...4..    .5.4..    .5.431    .54321
%    .....1    ...3.1    .4.3.1    54.3.1    54.310    543210
%NIL
%CL-USER(23): 
%\end{verbatim}

\section{Monodromy and Galois groups}\label{s:monodromy}

The monodromy group of a problem in enumerative geometry captures
information reflecting three aspects of algebraic geometry: geometry,
algebra, and arithmetic.  Informally, it is the symmetry group of the
set of solutions.  Three more precise interpretations are given below.
Historically, these groups were studied since the nineteenth century
\cite{jordan, dickson, weber}; modern interest probably dates from a
letter from Serre to Kleiman in the seventies (see the historical
discussion in the survey article \cite[p.~325]{kleiman}).  Their
modern foundations were laid by Harris \cite{harris}; among other
things, he showed that the monodromy group of a  problem
is equivalent to the Galois group of the equations defining it.

These groups are difficult to compute in general, and indeed they are
known for relatively few enumerative problems.  In this section, we
use the computation of explicit algebraic solutions to Schubert
problems (along with a criterion from \cite{VakilB}) to give a method
to compute many such groups explicitly (when they are ``full'', or as
large as possible), and to give an experimental method to compute
groups in other cases.

It is most interesting to exhibit cases where the Galois/monodromy
group is unexpectedly small.  Indeed, Harris writes of his calculations:
\begin{quote} the results represent an affirmation of
one understanding of the geometry underlying each of these problems, in
the following sense: in every case dealt with here, the actual
structure on the set of solutions of the enumerative problem as
determined by the Galois group of the problems, is readily described
in terms of standard algebrao-geometric constructions.  In particular,
in every case in which current theory had failed to discern any
intrinsic structure on the set of solutions --- it is proved here ---
there is in fact none. \cite[p.~687-8]{harris}
\end{quote}
We exhibit an example of a Schubert problem whose Galois/monodromy
group experimentally appears to be smaller than expected --- it is the
dihedral group $D_4 \subset S_4$.  This is the first example in which
current theory fails to discern intrinsic structure.  Examples of
``small'' Galois groups were given in \cite[Sect.\ 5]{VakilB}; but
there an explanation had already been given by Derksen.  Here,
however, we have a mystery: We do not understand geometrically why the
group is $D_4$.  (However, see the end of this section for a
conjectural answer.)

We now describe the three interpretations of the Galois/monodromy
group for a Schubert problem.
%corresponding to those $k$-planes in
%$n$-space ($K^n$) meeting $m$ general flags $F^i_\bullet$ ($1
%\leq i \leq m$) in given ways (corresponding to partitions $\al^1$,
%\dots, $\al^m$).  
The definition for a general problem in enumerative geometry is the
obvious generalization; see \cite{harris} for a precise definition,
and for the equivalence of \emph{(A)} and \emph{(B)}. See
\cite[Sect.\ 2.9]{VakilB} for more discussion.

{\em (A) Geometry.}  Begin with $m$ general flags; suppose
there are $N$ solutions to the Schubert problem (i.e.\ there are $N$
flags meeting our $m$ given flags in the specified manner).
%
%(i.e.\ $N$ different $k$-spaces satisfying the $m$ Schubert conditions). 
%
Move the $m$ flags around in such a way that no two of the solutions
ever come together, returning the $m$ flags to their starting positions, 
and follow the $N$ solutions.
The $N$ solutions are returned to their initial positions as a {\em
set}, but the individual $N$ solutions may be permuted.  What are the
possible permutations?  (See the applet
\verb+http://lamar.colostate.edu/~jachter/mono.html+ for an
illustration of this concept.)

{\em (B) Algebra.}  The $m$ flags are parameterized by $\flags^m$.
Define the ``solution space'' to be the subvariety of $\flags \times
\flags^m$ mapping to $\flags^m$, corresponding to those flags
satisfying the given Schubert conditions.  There is one irreducible
component $X$ of the solution space mapping dominantly to $\flags^m$;
the morphism has generic degree $N$.  The Galois/monodromy group is
the Galois group of the Galois closure of the corresponding extension
of function fields.  The irreducibility of $X$ implies that the Galois
group $G$ is a transitive subgroup of $S_N$.

{\em (C) Arithmetic.}  If the $m$ flags are defined over $\Q$, then
the smallest field of definition of a solution must have Galois group
that is a subgroup of the Galois/monodromy group $G$.  Moreover, for a
randomly chosen set of $m$ flags, the field of definition will have
Galois group precisely $G$ with positive probability (depending on the
particular problem).  The
equivalence of this
version with the previous two follows from \emph{(B)} by the Hilbert
irreducibility theorem, as $\flags^m$ is rational (\cite[Sect.\ 9.2]{lang},
see also \cite[Sect.\ 1.5]{serre} and \cite{cohen}).
We are grateful to M.~Nakamaye for discussions on this topic.

Given any enumerative problem with $N$ solutions, we see that the
Galois/ monodromy group is a subgroup of $S_N$; it is well-defined up
to conjugacy in $S_N$.  As the solution set should be expected to be as
symmetric as possible, one should expect it to be as large as
possible; it should be $S_N$ unless the set of solutions has some additional
geometric structure.
%%%%%%%%%%%%%%%%%%%%%%%%%%%%%%%%%%%%%%%%%%%%%%%%%%%%%%%%%%%%

For example, in \cite{harris}, Harris computed several
Galois/monodromy groups, and in each case they were the full symmetric
group, unless there was a previously known geometric reason why the
group was smaller.  The incidence relations of the $27$ lines on a
smooth cubic surface prevent the corresponding group from being
two-transitive.  There exist two of the $27$ lines that intersect,
and there exist another two that do not.  These incidence relations
can be used to show that the Galois/monodromy group must be contained
in the reflection group $W(E_6) \subset S_{27}$, 
e.g.\ \cite[Sects.\ 25,
26]{manin} or \cite[Prob.\ V.4.11]{hartshorne}; Harris shows that
equality holds \cite[III.3]{harris}.

Other examples can be computed based on permutation arrays.

\begin{corollary}\label{c:monodromy}
The explicit equations defining a Schubert problem in
Theorem~\ref{t:equations} can be used to determine the
Galois/monodromy group for the problem as well.
\end{corollary}

As a toy example, we see that the monodromy group for
Example~\ref{ex:2} is $S_{2}$, as there are two solutions to the
Schubert problem, and the only transitive subgroup of $S_2$ is $S_2$
itself.  Algebraically, this corresponds to the fact that the roots of
the irreducible quadratic $c^2+10c-20$ in example~\ref{ex:2} generate
a Galois extension of $\Q$ with Galois group $S_2$.
%
% ...as
%calculated using the command ``galois'' in Maple.  Note, the linear
%equations defining the problem do not impact the Galois/monodromy
%group for the problem.    <--- I don't know what this last sentence
% means!  -- Ravi

Unfortunately, the calculations of monodromy groups for flag varieties
becomes computationally infeasible as $n \rightarrow 10$ where the
number of solutions becomes larger.  Therefore, we have considered
related problems of computing Schubert problems for the Grassmannian
manifolds $G(k,n)$.  Here, $G(k,n)$ is the set of $k$-dimensional
planes in $\mathbb{C}^{n}$.  Schubert varieties are defined analogously by rank
conditions with respect to a fixed flag.  These varieties are indexed
by partitions $\lambda=(\lambda_{1}, \dotsc, \lambda_{k})$ where
$\lambda_{1}\geq \dotsb \geq \lambda_{k}\geq 0$.  The permutation
arrays work equally well for keeping track of the rank conditions for
intersecting Schubert varieties in the Grassmannian if we replace the
condition that a permutation array must have rank $n$ by requiring
rank $k$.

In the case of the Grassmannian, combinatorial criteria were given for
the Galois/monodromy group of a Schubert problem to be $A_N$ or $S_N$
in \cite{VakilB}.  Intersections on the Grassmannian manifold may be
interpreted as a special case of intersections on the flag manifold,
so our computational techniques apply.  We sketch the criteria here,
and refer the reader to \cite{VakilB} for explicit descriptions and
demonstrations.

\begin{criterion}
\label{criteria} \textsc{\textbf{Schubert Induction.}}  Given a
Schubert problem in the Grassmannian manifold, a choice of geometric
degenerations yields a directed rooted tree.  The edges are directed
away from the root.  Each vertex has out-degree between $0$ and $2$.
The portion of the tree connected to an outward-edge of a vertex is
called a {\em branch} of that vertex.  Let $N$ be the number of leaves
in the tree.
\begin{enumerate} 
\item[(i)]
Suppose each vertex with out-degree two satisfies either
(a) there are a different number of leaves on the two branches,
or (b) there is one leaf on each branch.  Then the Galois/monodromy
group of the Schubert problem is $A_N$ or $S_N$.
\item[(ii)]
Suppose each vertex with out-degree two has a branch with one leaf.
 Then the Galois/monodromy
group of the Schubert problem is $S_N$.
\item[(iii)]  Suppose that each vertex with out-degree two
satisfies (a) or (b) above, or (c) there  are $m \neq 6$ leaves
on each branch, and it is known that the corresponding Galois/monodromy
group is two-transitive.  Then the Galois/monodromy group is $A_N$ or
$S_N$.
\end{enumerate}
\end{criterion}

Part (i) is \cite[Thm.\ 5.2]{VakilB}, (ii) follows from the proof of
\cite[Thm.\ 5.2]{VakilB}, and (iii) is \cite[Thm.\ 5.10]{VakilB}.
Criterion (i) seems to apply ``almost always''.
Criterion (ii) applies rarely.  Criterion (iii) requires
additional information and is useful only in ad hoc circumstances.

\bigskip

%\subsection{Galois/monodromy groups from explicit algebraic solutions}
The method discussed in this paper of explicitly (algebraically)
solving Schubert problems gives two new means of computing Galois
groups.  The first, in combination with the Schubert induction rule,
is a straightforward means of proving that a Galois group is the full
symmetric group.  The second gives strong experimental evidence (but
no proof!) that a Galois group is smaller than expected.

\begin{criterion}\label{crit:2}
\textbf{\textsc{Criterion for Galois/monodromy group to be full.}}
If $m$ flags defined over $\Q$ are exhibited such that the solutions
are described in terms of the roots of an irreducible degree $N$
polynomial $p(x)$, and this polynomial has a discriminant that is not
a square, then by the arithmetic interpretation (C) above, the
Galois/monodromy group is not contained in $A_N$.  
\end{criterion}

Hence in
combination with the Schubert induction criterion (i), this gives a
criterion for a Galois/monodromy group to be the full symmetric group
$S_N$.

(In principle one could omit the Schubert induction criterion: if one
could exhibit a single Schubert problem defined over $\Q$ whose Galois
group was $S_N$, then the Galois/monodromy group would have to be
$S_N$ as well.  However, showing that a given degree $N$ polynomial
has Galois group $S_N$ is difficult; our discriminant criterion is immediate to
apply.)

The smallest Schubert problem where Criterion~\ref{criteria}(i)
applies but Criterion~\ref{criteria}(ii) does not is the intersection
of six copies of the Schubert variety indexed by the partition $(1)$
in $G(2,5)$ (and the dual problem in $G(3,5)$).  Geometrically, it
asks how many lines in $\proj^4$ meet six planes.  When the planes are
chosen generally, there are five solutions (i.e.\ five lines).  By
satisfying the first criterion we know the Galois/monodromy group is
``at least alternating'' i.e.\ either $A_N$ or $S_N$, but we don't
know that the group is $S_N$.  We randomly chose six planes defined
over $\Q$.  Maple found the five solutions, which were in terms of the
solutions of the quintic $101 z^5 - 554 z^4 + 887 z^3 -536 z^2 + 194 z
- 32$.  This quintic has non-square discriminant, so we conclude that
the Galois/monodromy group is $S_5$.  As other examples, the Schubert
problem $(2)^2 (1)^4$ in $G(2,6)$ has full Galois/monodromy group
$S_6$, the Schubert problem $(2) (1)^6$ in $G(2,6)$ has full
Galois/monodromy group $S_9$, and the Schubert problem $(2,2)(1)^5$ in
$G(3,6)$ has full Galois/monodromy group $S_6$.  We applied this to
many Schubert problems and found no examples satisfying
Criterion~\ref{criteria}(i) or (iii) that did not have full Galois
group $S_N$.

As an example of the limits of this method, solving the
Schubert problem $(1)^8$ in $G(2,6)$ is not computationally feasible
(it has $14$ solutions), so this is the smallest Schubert problem
whose Galois/monodromy group is unknown (although Criterion
~\ref{criteria}(i) applies, so the group is $A_{14}$ or $S_{14}$).

\begin{criterion}\label{crit:3}
\textbf{\textsc{Probabilistic evidence for smaller
Galois/monodromy groups.}} If for a fixed Schubert problem, a large
number of ``random'' choices of flags in $\Q^{n}$ always yield Galois
groups contained in a proper subgroup $G \subset S_N$, and the group
$G$ is achieved for some choice of Schubert conditions, this gives
strong evidence that the Galois/monodromy group is $G$.
\end{criterion}

This is of course not a proof --- we could be very unlucky in our
``random'' choices of conditions --- but it leaves little doubt.

As an example, consider the Schubert problem $(2,1,1) (3,1) (2,2)^2$
in $G(4,8)$.  There are four solutions to this Schubert problem.  When
random (rational) choices of the four conditions are taken, Maple
always (experimentally!) yields a solution in terms of $\sqrt{a + b
\sqrt{c}}$ where $a$, $b$, and $c$ are rational.  The Galois group of any such
algebraic number is contained in $D_4$: it is contained in $S_4$ as
$\sqrt{a+b \sqrt{c}}$ has at most $4$ Galois conjugates, and the
Galois closure may be obtained by a tower of quadratic extensions over
$\Q$.  Thus the Galois group is a $2$-subgroup of $S_4$ and hence
contained in a $2$-Sylow subgroup $D_4$.

%The only transitive subgroups of $S_4$ contained in $D_4$
%are (conjugate to) $C_2 \times C_2$, $C_4$, and $D_4$.

We found a specific choice of Schubert conditions for which the Galois
group of the Galois closure $K$ of $\Q(\sqrt{a+b \sqrt{c}})$ over $\Q$
was $D_4$.  (The numbers $a$, $b$, and $c$ are large and hence not
included here; the Galois group computation is routine.)  Thus we have
rigorously shown that the Galois group is at least $D_4$, hence $D_4$
or $S_4$.  We have strong experimental evidence that the group is
$D_4$.

\bigskip
\noindent 
{\bf Challenge:}  Prove that the Galois group of this Schubert problem
is $D_4$.
\bigskip

We conjecture that the geometry behind this example is as follows.
Given four general conditions, the four solutions may be labeled $V_1$, \dots, $V_4$ 
so that either (i)  $\dim (V_i \cap V_j) = 0$ if
$i \equiv j \pmod 2$ and 
$\dim (V_i \cap V_j) = 2$ otherwise, or 
(ii)  $\dim (V_i \cap V_j) = 2$ if
$i \equiv j \pmod 2$ and 
$\dim (V_i \cap V_j) = 0$ otherwise.
If (i) or (ii) holds then necessarily $G \neq S_4$, implying $G \cong D_4$.

This example (along with the examples of \cite[Sect.\ 5.12]{VakilB}) naturally
leads to the following question.  Suppose $V_1$, \dots, $V_N$
are the solutions to a Schubert problem (with generally chosen conditions).
Construct a rank table 
$$\left\{ \dim \left( \bigcap_{i \in I} V_i \right) 
\right\}_{I \subset \{ 1, \dots, n \}}.$$ In each known example,
the Galois/monodromy group is precisely the group of permutations
of $\{ 1, \dots, n \}$ preserving the rank table.

\bigskip
\noindent 
{\bf Question:}  Is this always true?

%\begin{remark}
{\em Remark.}
Schubert problems for the Grassmannian varieties were among the first
examples where the Galois/monodromy groups may be smaller than
expected.  The first example is due to H.  Derksen; the ``hidden
geometry'' behind the smaller Galois group is clearer from the point
of view of quiver theory.  Derksen's example, and other infinite
families of examples, are given in \cite[Sect.\ 5.13--5.15]{VakilB}.
%\end{remark}

\section{Acknowledgments}
We would like to thank Eric Babson and M.~Nakamaye for helpful discussions.
%\bibliography{my}

\end{document}